\documentclass[reqno]{amsart}
\usepackage{amssymb,amscd,verbatim, amsthm,graphics, color,latexsym,amsmath,multicol}
\usepackage[all,cmtip]{xy}
\usepackage{fancybox}
\usepackage{longtable}

\usepackage[top=3cm, bottom=3cm, left=3.2cm, right=3.2cm]{geometry}

\newtheorem{thm}[equation]{Theorem}

\newtheorem{lemma}[equation]{Lemma}
\newtheorem{Expectation}[equation]{Expectation}

\newtheorem{prop}[equation]{Proposition}
\newtheorem{cor}[equation]{Corollary}

\newtheorem{defn}[equation]{Definition}
\newtheorem{conj}[equation]{Conjecture}

\newtheorem*{thm*}{Theorem}
\newtheorem*{rmk*}{Remark}
\theoremstyle{remark}
\newtheorem{remark}[equation]{Remark}
\theoremstyle{remark}

\newtheorem{sketch}[equation]{Sketch of the proof}

\newcommand \fk[1]{{{\mathfrak #1}}}
\newcommand \C[1]{{\mathcal #1}}

\newcommand \wti[1]{{\widetilde {#1}}}

\newcommand\fg{\mathfrak g}

\newcommand \bC{{\mathbb C}}
\newcommand \bF{{\mathbb F}}

\newcommand \bR{{\mathbb R}}
\newcommand \RR{{\mathbb R}}
\newcommand \bZ{{\mathbb Z}}
\newcommand \ZZ{{\mathbb Z}}

\newcommand \bn{{\mathbf n}}

\newcommand \QQ{{\mathbb Q}}

\newcommand \Qlb{\bar{{\mathbb Q}}_\ell}
\newcommand \ra{{\rightarrow}}
\newcommand{\isoarrow}{{~\overset\sim\longrightarrow~}}

\newcommand\CN{{\C N}}
\newcommand\G{{\Gamma}}

\newcommand\Frob{\mathsf{Frob}}
\newcommand\Lie{\mathsf{Lie}}

\newcommand\CA{{\C A}}
\newcommand\CB{{\C B}}
\newcommand\CH{{\C H}}
\newcommand\CC{{\C C}}
\newcommand\CL{{\C L}}

\newcommand\CO{{\C O}}

\newcommand\CW{{\C W}}

\newcommand\CU{{\C U}}

\newcommand\ep{{\epsilon}}

\newcommand\om{{\omega}}
\newcommand\al{{\alpha}}

\newcommand \vg{\check{\fk g}}

\newcommand \vfh{\check{\fk h}}
\newcommand \ve{{\check e}}

\newcommand \vf{{\check f}}
\newcommand \vh{{\check h}}

\newcommand\ft{{\mathfrak t}}

\newcommand\Hom{\operatorname{Hom}}
\newcommand\End{\operatorname{End}}
\newcommand\Ind{\operatorname{Ind}}

\newcommand\ad{\mathbf{A}}
\newcommand\hG{G^\vee}
\newcommand\LGr{{}^LG}
\newcommand\LGad{{}^LG^{\rm ad}}

\newcommand\val{\mathsf{val}}

\newcommand\ba{\mathbf {a}}

\def\<{\langle} 
\def\>{\rangle}

\numberwithin{equation}{section}

\begin{document}

\author{Dan Ciubotaru}
\thanks{D. Ciubotaru was partially supported by EPSRC EP/V046713/1(2021).
}

\address{Dan Ciubotaru\\Mathematical Institute, Oxford University, Oxford OX2 6GG, UK
}
 \email{dan.ciubotaru@maths.ox.ac.uk}

\author{Michael Harris}
\thanks{M. Harris  was partially supported by NSF Grants DMS-2001369 and DMS-1952667, and by a Simons Foundation Fellowship, Award number  663678
}

\address{Michael Harris\\
Department of Mathematics, Columbia University, New York, NY  10027, USA}
 \email{harris@math.columbia.edu}
 
 \title[On the generalized Ramanujan and Arthur conjectures over function fields]{On the generalized Ramanujan and Arthur conjectures over function fields}

\maketitle

\begin{abstract} Let  $G$ be a simple group
over a global function field $K$, and let $\pi$ be a cuspidal automorphic representation of $G$.  Suppose 
$K$ has two places $u$ and $v$ (satisfying a mild restriction on the residue field cardinality), at which the group $G$ is quasi-split,
such that  $\pi_u$ is  tempered and $\pi_v$ is unramified and generic.   We prove that $\pi_w$ is tempered at all unramified places $K_w$ at which $G$ is unramified quasi-split.  

More generally, the set of unitary spherical representations is partitioned according to nilpotent conjugacy classes in the Lie algebra of $G$.   We show that if $\pi_v$ is in the set corresponding to the nilpotent  class $N$, and if $\pi_u$ satisfies an analogous hypothesis, then $\pi_w$ belongs to the same class $N$, where $w$ is as above.  These results are consistent with conjectures of Shahidi and Arthur.

The proofs use the Galois parametrization of cuspidal representations due to V. Lafforgue to relate the local Satake parameters of $\pi$ to Deligne's theory of Frobenius weights.  The main observation is that, in view of the classification of  unitary spherical representations, due to Barbasch and the first-named author, the theory of weights excludes almost all complementary series as possible local components of $\pi$.  This in turn determines the local Frobenius weights at all unramified places.  In order to apply this observation in practice we need a result of the second-named author with Gan and Sawin on the weights of discrete series representations.
\end{abstract}

\section{Introduction}   
\subsection{The generalized Ramanujan conjecture}
One of the most striking implications of Laurent Lafforgue's proof of the Langlands correspondence for $GL(n)$ over function fields is the Ramanujan Conjecture:  If $K$ is a global function field over a finite field and $\pi$ is a cuspidal automorphic representation of $GL(n,K)$ with unitary central character, then the local components of $\pi$ are tempered everywhere.  This is believed to be the case as well when $K$ is a number field, but Lafforgue's proof relies on Deligne's purity theorem (also known as the Weil Conjectures), whose applications over number fields are limited.

Vincent Lafforgue has used the geometry of varieties (or stacks) over finite fields to construct global Langlands parameters for 
cuspidal automorphic representations of $G(K)$ when $G$ is any connected reductive group.  In particular, Deligne's purity theorem is used at least implicitly in much of the theory on which his construction relies.   It is therefore natural to ask whether his results imply anything like the Ramanujan Conjecture for general $G$.  The first observation is that it has long been known that the most direct translation of the Ramanujan Conjecture fails for groups other than $GL(n)$ and its inner forms;  there are numerous constructions of cuspidal automorphic representations of classical groups that are not locally tempered almost everywhere.  The Arthur Conjectures, to which we return in the second part of this introduction, were designed in part to address this phenomenon and to characterize the discrepancy from temperedness systematically.   See \cite{Sa} for a review of the generalized Ramanujan Conjecture, primarily in the context of number fields.

Shahidi's $L$-packet conjecture provides an important insight into the question.  The conjecture, in its recent refinement \cite[Conjecture 1.5]{LS}, asserts that every Arthur
packet of a quasi-split $G$ is tempered if and only if it contains a locally generic member with respect to any Whittaker datum.  
In view of the Arthur Conjectures, in other words, genericity, in some form, is the natural condition for a version of the Ramanujan 
Conjecture valid for general $G$.   In his article \cite{Sh}, Shahidi showed\footnote{The article \cite{LS} of B. Liu and Shahidi treats more general Arthur packets, in connection with a conjecture due to Jiang.}  that 
\begin{thm*}\cite[Theorem 6.2]{Sh}  Assume $G$ is a quasi-split connected reductive group over the global field $K$.  The Arthur Conjectures, together with a generalization due to Clozel \cite{Cl}, imply that if $\pi$ is 
an everywhere locally generic cuspidal automorphic representation of $G$ with unitary central character, then the local components $\pi_v$ are tempered for all but finitely many $v$.
\end{thm*}
Shahidi states his theorem for number fields but it should be valid for function fields as well.\footnote{Laurent Clozel has pointed out that, in his review of Shahidi's article for Mathematical Reviews, Wee Teck Gan noted that the version of the Arthur Conjectures assumed in Shahidi's Theorem 6.2 includes the Ramanujan Conjecture for $GL(n)$.  For number fields this is quite far from being established, but for function fields this is due  to Laurent Lafforgue \cite[Th\'eor\`eme VI.10]{Laf02}.}

The purpose of this paper is to prove a version of this theorem when $K$ is a function field.   Our theorem is unconditional -- we assume neither the Arthur Conjectures nor
Clozel's generalization -- and we assume only that $\pi_v$ is generic at one place, and that $\pi_u$ is tempered at a second place.   This is consistent with the version of
the generalized Ramanujan Conjecture that predicts that {\it globally generic} cuspidal automorphic representations are everywhere tempered -- in \cite{Sa} Sarnak attributes this
version to the article \cite{HPS} of Howe and Piatetski-Shapiro.    At present we are required to make  a few additional technical assumptions.  Here is our statement.

\begin{thm*}[Theorem \ref{mainthm}, below]\label{t:main-intro}  Let $K$ be a function field over a finite field of characteristic $p$, and let
$G$ be a connected  absolutely simple group over $K$.  Let $\pi$ be a cuspidal automorphic representation of $G$.  Suppose 
\begin{enumerate}
\item The local component $\pi_v$ is generic and unramified at a place $v$ of $K$ such that $G(K_v)$ is an unramified quasi-split group
\footnote{Here and elsewhere, we use this expression as shorthand
for the condition that $G$, as an algebraic group over $K_v$, is unramified and quasi-split.}, other than the quasi-split form of a unitary group in $2n+1$ variables;
\item There is a place $u$ of $K$ such that $\pi_u$  is tempered;
\end{enumerate}
Then for every place $w$ of $K$ at which $\pi$ is unramified and $G(K_w)$ is unramified quasi-split, the local component $\pi_w$ is tempered.
\end{thm*}

\begin{rmk*}  Strictly speaking, $\pi$ is defined over a number field and the appropriate notion, both in the hypothesis and in the conclusion, is $\iota$-tempered for some embedding $\iota:  \bar{\QQ} \hookrightarrow \bC$.
For a discussion of the relevant notions of temperedness see Section \ref{iotatempered} and Theorem \ref{mainthm}.  
\end{rmk*}

The proof comes down to assembling three basic structural features of three distinct aspects of the theory of automorphic representations, together with a recent result of the second-named author with Gan and Sawin that serves to hold the disparate parts of the argument together.  The first element of the assembly is V. Lafforgue's assignment to $\pi$ of a compatible family of semisimple $\ell$-adic Langlands parameters 
$$\sigma = \sigma_{\ell,\pi}:   Gal(K^{sep}/K) \ra \LGr(\Qlb),$$
where $\LGr$ is the Langlands dual group.  The properties of $\sigma$ are recalled in \S \ref{sec_param}; the important point is that the restriction of $\sigma$ to the decomposition group at an unramified place $w$ belongs to the conjugacy class of the Satake parameter of $\pi_w$.

The second element is Deligne's theory of weights of $\ell$-adic representations of the Galois group of a function field over a finite field.  Together with L. Lafforgue's proof of the global Langlands correspondence for $GL(n,K)$, Deligne's theory implies that the composition of $\sigma$ with any $N$-dimensional linear representation  
$$\tau:  \LGr \ra GL(N,\Qlb),$$
a semisimple rank $N$ $\ell$-adic local system, decomposes uniquely as the direct sum of punctually pure local systems of  (integer) weights $w_\chi$, each shifted by a character $\chi$ of the Galois group of the constant field that can be thought of as the non-integral part of the summand.  A more precise statement is given in \eqref{chidecomp} below.   

If we knew that $\tau \circ \sigma$ were irreducible for some $\tau$,  it would be attached by the Langlands correspondence for $GL(N)$ to a cuspidal automorphic representation of $GL(N)$, and Laurent Lafforgue's results, combined with the semisimplicity of $G$, would imply that no non-trivial $\chi$ appears in the direct sum decomposition.  The genericity hypothesis (1) in the statement of Theorem \ref{mainthm} should imply that some $\tau\circ\sigma$ is irreducible, but this is not known in general.  Instead,
hypothesis (2) in the statement of Theorem \ref{mainthm} provides the anchor that guarantees that only integer weights occur in 
$\tau\circ\sigma$.   With this in hand, we can apply the classification of generic spherical unitary representations of the quasi-split group $G(K_v)$.   Since every irreducible spherical representation factors through the adjoint group, it is sufficient to consider the classification of the spherical unitary dual when $G(K_v)$ is adjoint. This is due to Barbasch and the first-named author and its application to the problem at hand takes up most of the first part of the article.  By analyzing $\tau\circ \sigma$ when $\tau$ is the adjoint representation, and one additional minuscule representation where needed, a case-by-case consideration of the classification shows that no generic complementary series representation is compatible with the requirement that the weights be integral.  This in turn implies that $\tau\circ\sigma$ is punctually pure of weight $0$, which implies Theorem \ref{mainthm}.   

 If $G(K_v)$ is the quasi-split form of a unitary group in $2n+1$ variables, the analysis of the generic spherical unitary locus  can only conclude that if the local component $\pi_v$ in Main Theorem is not  tempered, then the real part of its Satake parameter is precisely one explicitly determined parameter, $\frac 12\omega_n$, see Remark \ref{PSU}.

As a corollary of Theorem \ref{t:main-intro}, we obtain that every unramified local component $\pi_w$ is generic, provided $G(K_w)$ is {unramified quasi-}split   subject to the same restrictions as above.  

\subsection{The Arthur conjectures}

The Ramanujan Conjecture -- that any cuspidal automorphic representation of $GL(n)$ over a global field $K$ is everywhere locally tempered --  was proved by Laurent Lafforgue when $K$ is a function field and is known for certain classes of representations when $K$ is a number field.  
The analogous statement is false except when $G$ is $GL(n)$ or a closely related group.  Arthur's Conjectures \cite{A} provide a quantitative measure of the failure of square-integral automorphic representations, including cuspidal representations, to be locally tempered in terms of an extra datum, given as a nilpotent conjugacy class in the Lie algebra $\vg$ of the Langlands $L$-group $\LGr$ of $G$, or equivalently in terms of an $SL(2)$-triple in $\vg$.  We recall a simplified form of Arthur's Conjectures, as explained in Clozel's article \cite{Cl}; this will be sufficient for our purposes.   Let $K$ be a global function field over the finite field $k_1$, and let $W_K$ be the Weil group of $K$; in practice we will be able to replace $W_K$ by the global Galois group $\Gamma_K = Gal(K^{sep}/K)$.  A global {\it Arthur parameter} $\psi$ is then a homomorphism to the $L$-group:  
$$\psi:  W_K \times SL(2,C) \ra \LGr(C),$$
where $C$ is an algebraically closed field of characteristic zero.  The parameter $\psi$ is required to satisfy several natural properties; in particular the restriction of $\psi$ to $SL(2,C)$ is assumed to be algebraic.  Thus it differentiates to a Lie algebra homomorphism
\begin{equation}\label{LieN}  d\psi:  \mathfrak{sl}(2,C) = \Lie(SL(2))(C) ~\ra~ \Lie~\hG(C),
\end{equation}
where $\hG(C) \subset \LGr(C)$ is the Langlands dual group. We define
\begin{equation}\label{Npsi}  N_\psi = d\psi(\begin{pmatrix} 0 & 1 \\ 0 &0 \end{pmatrix}) \in \Lie~\hG(C).
\end{equation}
The nilpotent conjugacy class containing $N_\psi$ is the fundamental invariant of the Arthur parameter.

The parameters $\psi_1$ and $\psi_2$ are {\it equivalent} if one can be obtained from the other by conjugation by the Langlands dual group $\hG(C)$. 

Let $U \subset G(\ad_K)$ be an open compact subgroup, and let $S = S(U)$ be the finite set of places of $K$ at which $U$ does not contain a hyperspecial maximal compact subgroup.  Let $\CA_{disc}(G,U,C)$ be the space of discrete (square-integrable) automorphic forms on $G(K)\backslash G(\ad_K)/U$ with values in $C$; the square-integrability condition is the usual one when $C = \bC$, and we assume it makes sense over our chosen $C$ as well.  Let 
$$\CA_{disc}(G,C) = \varinjlim_U \CA_{disc}(G,U,C).$$
Let $\Psi(G,U)$ be the set of equivalence classes of Arthur parameters $\psi$ as above that are unramified outside $S$.  Arthur's Conjecture then asserts that
\begin{conj}[Arthur]\label{A1}
\begin{equation}
\CA_{disc}(G,U,C) = \oplus_{\psi \in \Psi(G,U)} \CA_{disc,U,\psi},
\end{equation}
\end{conj}
Here for each $\CA_{disc,U,\psi}$ there is a finite set of irreducible and necessarily unitary representations $\Pi_\psi = \{\pi_{\psi,i}\}$ of $G(\ad_K)$, $i \in I_\psi$, such that
$$\CA_{disc,\psi} = \oplus_i \pi_{\psi,i}^U.$$

This means in particular that for $w \notin S$, the local component $\pi_{\psi,i,w}$ is spherical as well as unitary.  The Conjecture also identifies the Satake parameters of $\pi_{\psi,i,w}$ in terms of the homomorphism
\begin{equation}\label{WDw} \psi_w:  WD_w \times SL(2,C) \ra \LGr(C)
\end{equation}
obtained (up to $\hG(C)$-conjugation) from $\psi$.  Here $WD_w$ is the Weil-Deligne group of $w$, and there is a dictionary that identifies parameters as in \ref{WDw} with homomorphisms
\begin{equation}\label{2SL2}  \psi_w:  W_w \times SL(2,C) \times SL(2,C) \ra \LGr(C),
\end{equation}
where $W_w$ is the local Weil group of the completion $K_w$ of $K$ at $w$.  We use the same notation $\psi_w$ for the two expressions \eqref{WDw} and \eqref{2SL2}; the context will make it clear which one we mean.   

Since $\psi_w$ is unramified, it factors (in the version \eqref{WDw}) through $W_w/I_w \times SL(2,C)$, where
$I_w$ is the inertia group.  Assume for the moment that $C = \bC$, and consider the map
\begin{equation}\label{j}  j:  W_w/I_w \ra SL(2,C); ~~j(\Frob_w) = \begin{pmatrix} |q_w|^{\frac{1}{2}} & 0 \\ 0 & |q_w|^{-\frac{1}{2}} \end{pmatrix}.
\end{equation}
Then
\begin{conj}[Arthur]\label{A2}   For any $\psi \in \Psi(G,U)$ and any $\pi = \pi_{\psi,i}$  
(a)  The Satake parameter $s(\pi_w)$ of $\pi_w$, for $w \notin S(U)$, has {\bf integral weight}, in the sense that, for any algebraic representation $\tau:  \LGr \ra GL(N)$, the image of the conjugacy class $s(\pi_w)$ under $\tau$ is a diagonalizable matrix in $GL(N,\CC)$, each of whose eigenvalues is a Weil $q$-number of some weight $w_{i,\tau} \in \ZZ$.

(b) Moreover, we have the identity
$$s(\pi_w) = \psi_w(\Frob_w,j(\Frob_w)),$$
up to conjugacy, where the argument $\Frob_w,j(\Frob_w)$ is viewed as an element of $W_w/I_w \times SL(2,\CC)$.

(c)  The restriction of $\psi_w$ to $W_w/I_w \times \{1\} \subset W_w/I_w \times SL(2,\CC)$ is {\it isobaric of weight $0$}; i.e., it lies in a maximal compact subgroup of $\LGr(\CC)$.

\end{conj}
This version of Arthur's Conjecture is a rewriting of Conjectures 3 and 2 and equation (2.5)  (in that order) of \cite[\S 2.3]{Cl}.  The three parts taken together imply that the weights $w_{i,\tau}$ in part (a) are completely determined by the restriction of $\psi_w$ to $\{1\} \times SL(2,\CC)$, and indeed by the nilpotent conjugacy class of $N_\psi$ , defined as in \eqref{Npsi}.  Since $N_\psi$ is defined globally, this means that, 

\begin{Expectation}\label{weights}  For any $\pi \in \CA_{disc,\psi}$ the set of weights $w_{i,\tau}$ is conjecturally the same for all unramified places $w$ and depends only on $N_\psi$.
\end{Expectation}

In the second part of this article we verify Expectation \ref{weights}, and a version of Conjecture \ref{A2}, for $\pi \subset \CA_{disc}(G,C)$ under certain local hypotheses.  
To this end, we apply the complete classification of unitary spherical representations of the groups $G(K_w)$, as developed in \cite{Ba, Ci-F4,Ci-E6,BC-E8,BC-gen,BC-ajm}.   
There it is proved that the Satake parameters of unitary spherical representations belong to a collection of subsets of the set of all Satake parameters, which we call {\it complementary series} and that the set of complementary series is in bijection with a set of nilpotent conjugacy classes.  The complementary series attached to the nilpotent class $\{0\}$ classifies generic unitary spherical representations.  

We use the same strategy as in the generic case to define Arthur parameters.  We  consider a discrete automorphic representation $\pi$ that satisfies a condition at one place that guarantees, again using \cite{De}, that all weights at unramified places are integral.  At a second place we assume that the local representation is unramified and belongs to the complementary series indexed by the nilpotent conjugacy class $N$.  Using the classification again, we then promote the semisimple parametrization of Lafforgue and Xue \cite{Laf18,X} to an Arthur parameter $\psi_\pi$ that satisfies all three parts of Conjecture \ref{A2} at all unramified places.

The work of Xue extends the parametrization of \cite{Laf18} to some irreducible components of $\CA_{disc}(G,C)$, but at present it is not known how to extend it in general.  The Conjectures \ref{discfinite} and \ref{paramdisc} cited in the statement of our main theorem asserts that this is possible; however, the argument applies to all discrete representations to which Xue's methods do apply.

\begin{defn}\label{centralpoint}
The spherical unitary dual of $G(K_v)$ is partitioned into complementary series $\mathsf{CS}(G(K_v),N)$, see  section \ref{s:CS}. Given a nilpotent conjugacy class $N$, let  $(e^\vee,h^\vee,f^\vee)$ with $e^\vee = N$ be an $SL(2)$-triple. The \emph{central point} of $\mathsf{CS}(G(K_v),N)$ is the $G^\vee$-conjugacy class of the Satake parameter $q^{h^\vee/2}$. The corresponding spherical representation is the $AZ$-dual of a tempered generic representation, where $AZ$ is the Aubert-Zelevinsky involution \cite{Aubert95}.  
\end{defn}

We can now state our main results in relation to the Arthur conjecture.  

\begin{thm}\label{mainthmArthur} Let $\pi \subset \CA_{disc}(G,\bar{\QQ})$ with $m_U(\pi) \neq 0$.   Assume Conjectures \ref{discfinite} and \ref{paramdisc}; alternatively, assume $\pi$ admits a parametrization as in \cite{Laf18}.   In particular, what follows is valid if $\pi$ is  cuspidal.  

Suppose that:
\begin{enumerate}
\item  There is a local place $v$ at which the group $G$ is split.
\item the local component $\pi_v$ at the chosen split place $v$ belongs to the complementary series for the nilpotent conjugacy class of $N$.
\item there is a place $u \in S$ such that $\pi_u$ is {\it geometric} in the sense that the Frobenius eigenvalues of its semisimple parameter are Weil $q$-numbers (see
Definition \ref{algebraic}).    
\end{enumerate}
Then for all $w \notin S$ at which $G$ is split, $\pi_w$ belongs to the complementary series for the nilpotent conjugacy class of $N$.  

Moreover, at all such places, including the place $v$, the real part of the Satake parameter either equals the central point for this complementary series or, for certain classes $N$,  belongs to one of the  \emph{extraneous points} listed in Theorem \ref{t:central}.
\end{thm}
The theorem is only stated when $G(K_v)$ is split; we have not worked out the classification of geometric points in complementary series for general quasi-split groups. Readers are referred to Shahidi's article \cite{Sh} for indications of the compatibility of our results with the Arthur Conjectures.   Residual points can exist in a minority of complementary series.  The complete list of $N$ that admit extraneous points is given in Theorem \ref{t:central}.  The Arthur conjectures imply that such representations do not contribute to the discrete automorphic spectrum.  

Finally, we can eliminate the extraneous points if we make a stronger hypothesis on $\pi_u$.

\begin{thm*}[Theorem \ref{AZs}]  Under the hypotheses of Theorem \ref{mainthmArthur}, suppose in addition that $\pi_u = AZ(\sigma)$, where $\sigma$ is a discrete series representation with Weil-Deligne parameter $(\varphi_\sigma,N)$ (see \cite[Theorem 1.2]{GHS}) and $AZ$ denotes the Aubert-Zelevinsky involution, where $N$ determines the complementary series for the chosen  unramified place $v$.
Then at every split unramified place $w$, the real part of the Satake parameter of $\pi_w$ equals the central point for the complementary series for the nilpotent conjugacy class of $N$.  

In particular, no unramified local parameter of $\pi$ is an extraneous point.
\end{thm*}

\section*{Acknowledgments}  The authors thank Anne-Marie Aubert, Laurent Clozel, Johan de Jong, H\'el\`ene Esnault, Dennis Gaitsgory, Mark Kisin, Colette Moeglin, Peter Sarnak, Freydoon Shahidi,  Jack Thorne, and Xinwen Zhu for helpful clarifications and for help with references; in particular, we thank Zhu for pointing out his paper with D. Xu that allowed us to remove a slight restriction from the statement of our main result.  We also thank Martin Olsson for explaining the current state of the literature on weights in the $\ell$-adic cohomology of stacks, and for providing the sketch of the proof of the purity result used to provide an alternative approach to Theorem  \ref{mainthm}. 

Finally, we thank the anonymous referees whose comments encouraged us to include a sketch of the proof of a statement on the invariance of discrete series representations under Galois conjugation, and Peter Sarnak for asking the question that led to our results related to the Arthur conjectures.

\section{Nilpotent orbits and weight filtrations}

\subsection{Review of nilpotent orbits in semisimple Lie algebras}

Let $\fg^\vee$ be a semisimple Lie algebra over an algebraically closed field $C$ of characteristic zero.  (In practice, $C$ will be either $\bC$ or $\Qlb$.)  We let $G^\vee$ be a semisimple algebraic group over the field $C$, so that $\fg^\vee = \Lie(G^\vee)$. Let $\CN^\vee$ denote the nilpotent cone of $\fg^\vee$, the variety of nilpotent elements of $\fg^\vee$. It is well known that the group $G^\vee$ acts on $\CN^\vee$ with finitely many orbits, see for example \cite[Theorem 3.5.4]{CM93}. 

An $SL(2)$-triple in $\fg^\vee$ consists of elements $e^\vee, f^\vee, h^\vee \in \fg^\vee$ that satisfy the usual relations of the Lie algebra of $SL(2)$:
\begin{equation}\label{sl2}
[e^\vee,f^\vee] = h^\vee, [h^\vee,e^\vee] = 2e^\vee;  [h^\vee,f^\vee] = -2f^\vee.
\end{equation}
Notice that $e^\vee$ is necessarily nilpotent. Conversely, by a theorem of Jacobson and Morozov \cite[Theorem 3.3.1]{CM93}, if $e^\vee\in \CN^\vee$, there exists an $SL(2)$-triple $e^\vee, f^\vee, h^\vee \in \fg^\vee$. Moreover, Kostant proved that this gives in fact a one-to-one correspondence between $\hG$-orbits in $\CN^\vee$ and $G^\vee$-conjugacy classes of $SL(2)$-triples, \cite[Theorem 3.4.10]{CM93}.

\subsection{The weight filtration attached to an $SL(2)$-triple}

Let $Rep(\fg^\vee)$ denote the symmetric monoidal  category of pairs $(V,\rho)$ where $V$ is a finite-dimensional $C$-vector space
and $\rho:  \fg^\vee \ra \End(V)$ is a homomorphism of Lie algebras. Let $Fil_C$ denote the symmetric monoidal category of pairs $(V,Fil_\bullet(V))$ where  and $\dots Fil_i(V) \subset Fil_{i+1}(V) \dots$ is an increasing $\ZZ$-filtration. 
\begin{defn}\label{d:sl2-fil}  A weight filtration on $Rep(\vg)$ is a symmetric monoidal functor 
\begin{equation}
\mathbb{W}:  Rep(\fg^\vee) \ra Fil_C
\end{equation}
such that, for any $(V,\rho) \in Rep(\fg^\vee)$, $\mathbb{W}((V,\rho)) = (V,Fil_\bullet(V))$.
\end{defn}

A weight filtration is of {\it $SL(2)$-type} if there is an $SL(2)$-triple $e^\vee, f^\vee, h^\vee \in \fg^\vee$ such that, for every $(V,\rho)$ and all $i \in \ZZ$,
\begin{enumerate}
\item $e^\vee(Fil_i(V)) \subset Fil_{i+2}(V)$;
\item $h^\vee(Fil_i(V)) \subset Fil_i(V)$ and $h^\vee$ has eigenvalue $i$ on $gr_i(V)$.
\end{enumerate}

Any $SL(2)$-triple in $\fg^\vee$ defines a weight filtration on $Rep(\fg^\vee)$, necessarily of $SL(2)$-type.

\begin{lemma}\label{l:sl2-fil}  Let $(e^\vee, f^\vee, h^\vee)$, $(e'^\vee, f'^\vee, h'^\vee)$ be two $SL(2)$-triples in $\fg^\vee$, and let $(\rho,V) \mapsto (V,Fil_\bullet(V))$, $(\rho,V) \mapsto (V,Fil'_\bullet(V))$ be the corresponding weight filtrations on $Rep(\fg^\vee)$.  Suppose for all $(\rho,V)$, and for all $i \in \ZZ$, 
$$\dim Fil_i(V) = \dim Fil'_i(V).$$
Then $(e^\vee, f^\vee, h^\vee)$, $(e'^\vee, f'^\vee, h'^\vee)$ are conjugate under the adjoint action of the semisimple group $G^\vee$.
\end{lemma}

\begin{proof}
This is a particular case of a more general statement, see Proposition \ref{p:weight-pattern}. In Appendix, we also give an explicit proof of this result using the defining representation for the classical simple Lie algebras and the adjoint representation for the exceptional Lie algebras.
\end{proof}

\subsection{Frobenius weight filtrations}

Let $F = k((t))$ be a local field of positive characteristic $p$, with $k$ a finite field of order $p^r$ for some $r$.  Let $\overline F$ be an algebraic closure of $F$ and $F^{sep}$ the separable closure of $F$ in $\overline F$. Let $\G = \G_F = Gal(F^{sep}/F)$ and let 
${}^LG = \hG \rtimes \G$ be the $L$-group of a reductive group $G/F$, see \cite[\S1]{Kot} for example. More precisely, fix a maximal torus $T$, a Borel subgroup $B\subset T$ of $G(\overline F)$, and let $X^*(T)$ (resp. $X_*(T)$) be the character (resp cocharacter) group of $T$. Let $\Delta^*(T,B)$ (resp $\Delta_*(T,B)$) be the set of simple $B$-positive roots (resp. coroots) of $T$. The quadruple $\Psi_0(G(\overline F))=(X^*(T), X_*(T),\Delta^*(T,B),\Delta_*(T,B))$ is a based root datum of $G(\overline F)$. The Galois group $\Gamma$ acts on $\Psi_0(G(\overline F))$.

The L-group ${}^LG$ is a connected reductive group $G^\vee$ together with an action of $\Gamma$ on $G^\vee$ such that 
\begin{enumerate}
\item[(i)] there exists a $\Gamma$-isomorphism between $\Psi_0(G^\vee)$ and the dual of $\Psi_0(G(\overline F))$, and
\item[(ii)] there exists a {\it splitting} $(T^\vee,B^\vee, \{X_{\alpha^\vee}\})$ of $G^\vee$ fixed by $\Gamma$, where $T^\vee$ is a maximal torus of $G^\vee$, $B^\vee\supset T^\vee$ is a Borel subgroup and $X_{\alpha^\vee}\neq 0$ is an element in the root space $\mathfrak g^\vee_{\alpha^\vee}$, $\alpha^\vee\in \Delta^*(T^\vee,B^\vee)\cong\Delta_*(T,B)$.
\end{enumerate}

Let
$$\varphi:  \G \ra {}^LG$$
be an $L$-homomorphism.  For any finite-dimensional representation $\tau:  {}^LG \ra GL(N)$, let $(\rho(\varphi,\tau),L(\varphi,\tau))$ denote the representation $\tau\circ \varphi$.  
\begin{defn}\label{algebraic}  Let $\Frob \in \G$ be any Frobenius element.  We say $\varphi$ is {\bf geometric} if for one (equivalently for any) faithful representation $\tau$ as above, the eigenvalues of $\rho(\varphi,\tau)(\Frob)$ are Weil $q$ numbers of integer weight.  In other words, if $\lambda$ is any eigenvalue of $\rho(\varphi,\tau)(\Frob)$, then 
$\lambda$ is an algebraic number and there is an integer $w$ such that, for any embedding $\iota: \QQ(\lambda) \hookrightarrow \CC$, 
$$|\iota(\lambda)| = q^{\frac{w}{2}}.$$
\end{defn}

The definition is identical if $\CC$ is replaced by any algebraically closed field $C$ of characteristic zero; we will use it without comment when $C = \Qlb$.


The following definition will play an important role in the proofs.

\begin{defn}\label{mwf} Let  $w$ be an integer.  A  \emph{pure Weil-Deligne representation of weight $w$} of $K_u$ is a pair $(\varphi,N)$, where 
$$\varphi:  W_{K_u} \ra GL(V)$$
is a representation of the Weil group of the local field $K_u$ on an $m$-dimensional $\Qlb$-vector space $V$, and $N: V \ra V$ is a nilpotent operator, such that
the pair $(\varphi,N)$ satisfies the usual properties of a Weil-Deligne representation, and such that, for any choice of Frobenius element $\Frob_u \in W_{K_u}$ we have
\begin{itemize}
\item[(i)]  The eigenvalues of $\varphi(\Frob_u)$ are all $q$-numbers of integer weight.
\item[(ii)]  The subspace $W_aV \subset V$ of eigenvectors for $\varphi(\Frob_u)$ with eigenvalues of weight $\leq a$ is invariant
under $(\varphi,N)$;
\item[(iii)]  Letting $gr_aV = W_aV/W_{a-1}V$ the map
$$N^i:  gr_{w-i}V \ra gr_{w+i}V$$
is an isomorphism for all $i$.
\end{itemize}
\end{defn}

\subsection{Weight filtrations and local systems on curves}

In what follows, $X$ is a complete curve over $Spec(k)$, $S$ a non-empty finite set of places of the function field $K = k(X)$, $|S| \subset X$ the support of $S$, $U = X \setminus |S|$.
The following is a direct consequence of Deligne's results on weights of local systems in \cite{De}.  The final statement is Deligne's theorem on the purity of the monodromy weight filtration.

\begin{prop}\label{globalweights}  Let $\CL$ be an $\ell$-adic local system on $U$, $u \in S$, $K_u$ the completion of $K$ at $u$, and let $\Gamma_u = Gal(\bar{K}_u/K_u)$.  Choose a Frobenius element $\Frob_u \in \Gamma_u$.
Let $\rho_u:  \Gamma_u \ra Aut(\CL_u)$ denote the monodromy representation at $u$ attached to $\CL$.   Suppose the representation $\rho_u$ is geometric, in the sense of Definition \ref{algebraic} above.  In other words, suppose every eigenvalue $\alpha_j$ of $\rho_u(\Frob_u)$ is an algebraic number such that, for every  embedding $\beta:  \QQ(\alpha_j) \hookrightarrow \CC$,
$|\beta(\alpha_j)| = q^{\frac{w_j}{2}}$ for some integer $w_j = w(\alpha_j)$.  

Then for every $w \notin S$, the (unramified) monodromy representation $\rho_w: \Gamma_w \ra Aut(\CL_w)$ is also geometric.  Let $\gamma_j$ be the set of eigenvalues of $\rho_w(\Frob_w)$, for a Frobenius element $\Frob_w \in \Gamma_w$; let $\CW(\CL) \subset \ZZ$ be the set of Frobenius weights $w_i$ that occur in $\CL_w$, and, for each $i \in \CW(\CL)$, let $n(i)$ be the number of $\gamma_j$ such that $w(\gamma_j) = i$.

Moreover, $\CL$ has an increasing filtration $W_\bullet(\CL)$ such that, for each $i \in \ZZ$, 
$$gr_i(\CL) := W_i(\CL)/W_{i-1}(\CL)$$
is punctually pure of weight $i$ and $\dim gr_i(\CL) = n(i)$.

Finally, with $u \in S$ as above, let $\rho_{u,i}$ denote the monodromy representation of $\Gamma_u$ on the fiber $gr_i(\CL)_u)$.  Then $\rho_{u,i}$ is a pure Weil-Deligne representation of weight $i$.
\end{prop}


\section{Parameters of cuspidal representations and Weil numbers}\label{sec_param}

In what follows $k_1$ is a finite field of order $q_1 = p^f$ for some prime $p$, $X$ is an algebraic curve over $k_1$ with function field
$K = k_1(X)$, $v$ is a chosen place of $K$ with residue field $k$ of order $q$, and $G$ is a semisimple group over $K$ that is {unramified quasi-}split at $v$.  {For
reasons to be explained in Remark \ref{PSU}, we assume that $G(K_v)$ is not isomorphic to the quasi-split but non-split group $PSU(2n+1)$. }

We will fix
an open compact subgroup $U \subset G(\ad_K)$ that contains a hyperspecial maximal compact subgroup $U_v \subset G(K_v)$.   Let $X_U$ denote the double coset space $G(K)\backslash G(\ad_K)/U$.
The space of automorphic forms  of level $U$ with coefficients in a (commutative) ring $E$ is denoted
$$\CA(G,U,E) := \{f:  X_U \ra E \},$$
where all $f$ are assumed continuous with respect to the discrete topology on $E$.  
 
Assume $p$ is invertible in $E$ and let $\CA_0(G,U,E) \subset \CA(G,U,E)$ denote the submodule of cusp forms.  This is defined by the vanishing of
constant terms, which are defined by integrals with respect to a measure with values in $\ZZ[\frac{1}{p}]$.   
It is clear that if $E$ is a subring of the algebraically
closed field $C$
then the inclusion of $E$ in $C$ defines an isomorphism
$$\CA_0(G,U,E)\otimes C \isoarrow \CA_0(G,U,C).$$
We let 
$$\CA_0(G,E) = \varinjlim_{U} \CA_0(G,U,E),$$
the colimit taken with respect to inclusion.  
 
 \begin{prop}\label{cuspfinite}  Suppose $E$ is a noetherian $\ZZ[\frac{1}{p}]$ algebra with a fixed embedding in the algebraically
 closed field $C$.
 Then $\CA_0(G,U,E)$ is a finite $E$-module and the functions in $\CA_0(G,U,E)$ are compactly supported in $X_U$.
 Moreover, $\CA_0(G,U,E)\otimes C \isoarrow \CA_0(G,U,C).$
  \end{prop}
  
  \begin{proof}  All but the final statement is contained in \cite[Proposition 8.2]{BHKT} (which uses the assumption that $G$ is semisimple); the last statement is obvious.
  \end{proof}
  
Henceforward $E$ will be a number field and $C$ will be either $\bC$ or $\Qlb$.  We fix an algebraic closure
  $\bar{\QQ}$ of $\QQ$ containing $E$ and choose embeddings
  $$\iota:  \bar{\QQ} \hookrightarrow \bC; ~~~ \iota_\ell:  \bar{\QQ} \hookrightarrow \Qlb.$$

\subsection{Rationality of cuspidal representations}\label{rcr}  Let $S = S(U)$ be the finite set of places $w$ of $K$ such that $U_w := U \cap G(K_w)$ is not a hyperspecial
maximal compact subgroup; thus $v \notin S$.   We define the Hecke algebra 
$$\CH^S(G,U) = \otimes'_{w \notin S} \CH_w$$
where $\CH_w = \CH(G(K_w),U_w)$ is the algebra of $\QQ$-valued compactly supported functions on $G(K_w)$, biinvariant under $U_w$, with multiplication given
by convolution; the restricted tensor product $\otimes'$ is taken with respect to the unit elements $\mathbf{1}_{U_w} \in \CH_w$.   This is a commutative algebra and for any
sufficiently large number field $E$ as in Proposition \ref{cuspfinite} (depending on $U$) the action of $\CH^S(G,U)\otimes E$ on $\CA_0(G,U,E)$ is semisimple and decomposes $\CA_0(G,U,E)$ as a finite sum of
eigenspaces.  (Here again we are using our assumption that $G$ is semisimple.)  It follows in particular that 

\begin{lemma}\label{algHecke}  The eigenvalues of $\CH^S(G,U)$ on $\CA_0(G,U,C)$ are algebraic numbers for any algebraically closed $C$ of characteristic zero.
\end{lemma}

Similarly, the action of $G(\ad_K)$ on cusp forms preserves the rational subspace $\CA_0(G,\QQ)$.  So any  irreducible cuspidal automorphic representation has a model
over a number field.  It follows that 
$$\CA_0(G,\bar{\QQ}) \isoarrow \bigoplus_{\pi}  m(\pi) \pi;   \CA_0(G,U,\bar{\QQ}) \isoarrow \bigoplus_{\pi} m_U(\pi) \pi^U$$
where $\pi$ runs over the irreducible $\bar{\QQ}$-representations of $G(\ad_K)$ and $m(\pi) \geq 0$ is an integer multiplicity, and we make the convention
that $m_U(\pi) = 0$ if $\pi^U = 0$.  

Fix a prime $\ell \neq p$ and a level subgroup $U$.  Let $\CH^S(G,U,\Qlb) = \CH^S(G,U)\otimes \Qlb$.   In \cite{Laf18}
Vincent Lafforgue defines a commutative algebra $\CB(G,U,\Qlb) \supset \CH^S(G,U,\Qlb)$ of {\it excursion operators} that acts on $\CA_0(G,U,\Qlb)$ and defines a
decomposition
\begin{equation}\label{param}
\CA_0(G,U,\Qlb) \isoarrow \bigoplus_{\sigma} \CA_0(G,U,\Qlb)_{\lambda(\sigma)},
\end{equation}
where $\lambda(\sigma)$ runs over characters of $\CB(G,U,\Qlb)$ that occur non-trivially in $\CA_0(G,U,\Qlb)$, and $\sigma$ designates a semisimple Langlands parameter
\begin{equation}\label{sigma}
\sigma = \sigma_\ell:  Gal(K^{sep}/K) \ra \LGr(\Qlb).
\end{equation}
It is not known in general whether or not this decomposition is defined over $\bar{\QQ}$, but it is known for $GL(n)$ (when an appropriate central character is fixed).  However, 
\begin{lemma}\label{excHecke}  The restriction of the action of $\CB(G,U,\Qlb)$ on $\CA_0(G,U,\Qlb)$ to the subalgebra $\CH^S(G,U,\Qlb)$ coincides with the $\Qlb$-linear
extension of the natural action of $\CH^S(G,U,\bar{\QQ})$.
\end{lemma}
This is a special case of Proposition 0.16 (Proposition 6.2) of \cite{Laf18}; it is  one of the key steps in Lafforgue's construction.

Let $\Phi(U)$ denote the set of $\sigma$ such that $\CA_0(G,U,\Qlb)_{\lambda(\sigma)} \neq 0$.    For $\sigma \in \Phi(U)$, we let $\lambda_\CH(\sigma)$ denote
the restriction of $\lambda(\sigma)$ to $\CH^S(G,U,\bar{\QQ})$; for $w \notin S(U)$ we let  $\lambda_w(\sigma)$ denote the restriction of $\lambda(\sigma)$ to $\CH_w$.  
We let $\LGr^{ss}(\Qlb)$ denote the set of semisimple elements of $\LGr(\Qlb)$, $[\LGr^{ss}(\Qlb)]$ the set of semisimple conjugacy classes in $\LGr(\Qlb)$.
For $\sigma \in \Phi(U)$ and $w \notin S(U)$, let 
$$\alpha_w(\sigma) = [\sigma(\Frob_w)] \in [\LGr^{ss}(\Qlb)],$$
where the brackets around $\sigma(\Frob_w)$ denote the conjugacy class of the image under $\sigma$ of a choice of Frobenius element at $w$.

The key fact about Lafforgue's parametrization that we need is
\begin{thm}[Lafforgue]\label{SeqT}  Fix $\sigma \in \Phi(U)$.  Let
$s_{w,\ell}(\sigma) \in [\LGr^{ss}(\Qlb)]$ be the Satake parameter corresponding to $\lambda_w(\sigma)$.  Then 
$$s_{w,\ell}(\sigma) = \iota_\ell(s_w(\sigma))$$
for an $s_w(\sigma) \in [\LGr^{ss}(\bar{\QQ})]$
and
$$\alpha_w(\sigma) = s_{w,\ell}(\sigma).$$
\end{thm}

The first claim is a consequence of Lemma \ref{algHecke}; the second is contained in Th\'eor\`eme 0.1 (Th\'eor\`eme 11.11) of \cite{Laf18}. 

On the other hand, $\lambda_v(\sigma)$, with our fixed unramified place $v$, is the action of $\CH_v$ on $\pi^U$ for some $\pi \subset \CA_0(G,\bar{\QQ})$ with $m_U(\pi) \neq 0$.
Since $\pi$ is cuspidal, it is unitary, and therefore $\iota(s_v(\sigma)) \in [\LGr^{ss}(\bC)]$ is the Satake parameter of a unitary spherical 
representation $\pi_v$ of $G(K_v)$.    Let $T^\vee$ be a maximal torus of the Langlands dual group $\hG$ of $G$.  We may take the Satake parameter in $[\LGr^{ss}(\bC)]$ to be of the form $(s_v(\sigma)^0,\Frob_w)$ where 
$s_v(\sigma)^0 \in [T^{\vee}(\bar{\QQ})]$, the set of algebraic semisimple conjugacy classes in $\hG(\bar{\QQ})$; the brackets denote conjugacy classes under the Weyl group of $\hG$.\footnote{When $G$ is not split the parameter may only be determined up to the action of the Weil group, but none of what follows depends on this action.}
It follows from Theorem \ref{SeqT} that we thus have an algebraic semisimple conjugacy class $s_v(\sigma)^0$ in $T^\vee$ whose image under a complex embedding corresponds to a unitary spherical representation, while the image of $s_v(\sigma)$ under an $\ell$-adic embedding corresponds to the conjugacy class of the Frobenius at $v$ of
a semisimple Galois parameter with values in $\LGr(\Qlb)$.  The interaction between these properties of the  conjugacy class 
$s_v(\sigma)$ is the subject of our main theorems.

\begin{defn}\label{belongs}  Let $\pi \subset \CA_0(G,\bar{\QQ})$ with $m_U(\pi) \neq 0$.  We say $\pi$ {\bf belongs to} the parameter $\sigma$ if the restriction
$\lambda_\CH(\sigma)$ coincides with the action of $\CH^S(G,U,\bar{\QQ})$ on $\pi^U$.  
\end{defn}

\subsection{Local parametrizations}

We recall the local version of Lafforgue's parametrization of cuspidal automorphic representations.    Let $F = k((t))$ as above and
let $\Pi(G/F)$ be the set of (equivalence classes of) irreducible representations of $G(F)$ over the field $C$.  Let $\Phi(G/F)$ be the set of  (equivalence classes of) {\it semisimple } $L$-homomorphisms $W_F \ra {}^LG(C)$.    We take $C = \Qlb$ for some prime $\ell$.  

\begin{thm}[Genestier-Lafforgue, Fargues-Scholze, Li Huerta]\label{sslocal} 
There is a  semisimple local parametrization
$$\CL_F:  \Pi(G/F) \ra \Phi(G/F)$$ 
with the following properties.   
\begin{enumerate}
\item   Suppose $F = K_w$ where $K = k_1(X)$ and $w$ is any place of $K$.  Let $\Pi$ be a cuspidal automorphic representation of $G(\ad_K)$.  Then
$$\CL_F(\Pi_w) := [\CL(\Pi) ~|_{W_{F}}]^{ss}$$
depends only on $F = K_w$ and $\Pi_w$ (not on the globalizations $K$ and $\Pi$).     

\item $\CL_F$ is compatible with parabolic induction in the obvious sense.   
\item  In particular, if $\CL_F(\pi)$ is irreducible then $\pi$ is supercuspidal.

\end{enumerate}
\end{thm}

Li Huerta \cite{LH} has proved that the parametrizations of Genestier-Lafforgue \cite{GLa} and Fargues-Scholze \cite{FS} coincide for $F$ of positive characteristic.  

\begin{defn}\label{d:geometricrep}
We define a representation $\pi \in \Pi(G/F)$ to be \emph{geometric} if it has a model $\pi_E$ over a number field $\iota:  E \hookrightarrow \Qlb$:
$$\pi_E \otimes_{E,\iota} \Qlb \isoarrow \pi$$
such that the parameter $\CL_F(\pi):  W_F \ra {}^LG(\bar{\Qlb})$
is geometric in the sense of Definition \ref{algebraic}.  
\end{defn}

\begin{prop}\label{geometricpi}  Let $\pi \in \Pi(G/F)$ be given.  Then

\begin{enumerate}
\item  If $\pi$ is geometric, for any embedding $\iota':  E \ra \Qlb$, the representation $\pi_E \otimes_{E,\iota'}\Qlb$ is again geometric.  
\item  If $G$ is a torus then $\pi$ is geometric if and only if $\pi$ takes values in the group of Weil $q$-numbers.
\item  Geometricity is preserved under normalized parabolic induction, Jacquet modules,  direct products, and twisting by geometric characters.
\item  For any reductive $G$, if $\pi$ is geometric then so is its central character.
\item  If $\pi$ is discrete series with geometric central character then $\pi$ is geometric.  This is the case in particular for every discrete series when $G$ is semisimple.
\end{enumerate}

\end{prop}
\begin{proof}  The first point is a consequence of the compatibility of the Genestier-Lafforgue and Fargues-Scholze parametrizations with automorphisms of $\Qlb$.  The second can be reduced to the case where $G = \mathbb{G}_m$, in which case it follows from the compatibility of the parametrizations with class field theory.  The first three parts of point 3 are obvious, and the last follows from point 2.   Point 4 is automatic once one remembers that, under Langlands functoriality, restricting $\pi$ to the center $Z_G$ of $G$ corresponds to restricting the set of $\tau$ in the definition of geometricity to those $\tau$ factoring through the $L$-group of $Z_G$, which is a quotient of ${}^LG$.  

Finally point 5 follows from Theorem \ref{iotads} below.  

\end{proof}

\section{Some representation theory}

\subsection{Tempered and $\iota$-tempered representations}\label{iotatempered}

Let $F = k((t))$; in other words $F = K_v$ in the previous notation, and let $\pi$ be an irreducible admissible representation of $G(F)$ with coefficients in $\bC$.  Then $\pi$ is tempered
if and only if there is a parabolic subgroup $P \subset G(F)$ with Levi component $M$ and a discrete series representation $\rho$ of $M$, with unitary central character,
such that $\pi$ is an irreducible constituent of the (completely reducible) parabolically unitarily induced representation $\Ind_P^{G(F)} \rho$.  

\begin{defn}\label{tempitemp} Let $\pi$ be an irreducible admissible representation of $G(F)$ with coefficients in $\bar{\QQ}$.  Then we define $\pi$ to be {\it $\iota$-tempered} if $\pi \otimes_{\bar{\QQ},\iota} \bC$ is tempered  in the above sense.  
We say $\pi$
is {\it tempered} if it is $\iota'$-tempered for every embedding $\iota':   \bar{\QQ} \hookrightarrow \bC$.
\end{defn}

Suppose $\pi$ as above is $\iota$-tempered, and let $Z_M$ be the center of $M$.  Then the central character 
$$\xi_\rho:  Z_M \ra \bar{\QQ}^\times,$$
has the property that $\iota\circ \xi_\rho:  Z_M \ra \bC^\times$ is unitary.  There is no reason a priori for this to remain true if $\iota$ is replaced by a second embedding
$\iota':  \bar{\QQ} \hookrightarrow \bC$.    However, the following theorem is proved in \cite{GHS}:  

\begin{thm}\label{iotads}  Let $\pi$ be a complex-valued discrete-series representation of $G(F)$.  Then $\pi$ has a model $\pi_E$ over a number field $E \subset \bar{\QQ} \subset \bC$.  Moreover, let $\sigma \in Gal(\bar{\QQ}/\QQ)$, and let ${}^\sigma(\pi) = \sigma(\pi_E)\otimes_{\sigma(E)}\bC$.  Then ${}^\sigma(\pi)$ is in the discrete series.
\end{thm}

This is in fact a special case of the following result:

\begin{thm}\label{iotat}  Let $\Pi$ be a cuspidal automorphic representation of $G(\ad_K)$.   Then $\Pi$ has a model over a number field.  Suppose $\Pi_v$ is $\iota$-tempered for some place $v$ of $K$.  
Then $\Pi_v$ is in fact tempered.

\end{thm}

Thus we are entitled to omit the $\iota$- when referring to tempered components of cuspidal automorphic representations.

Complete proofs of the two theorems will appear elsewhere.  For the reader's convenience we supply a sketch of the proofs.  
\begin{enumerate}
\item We have seen in \S \ref{rcr} that a cuspidal $\Pi$ has a model over a number field.  It is proved in Beuzart-Plessis's appendix to \cite{GHS} that any discrete series representation occurs as a local component of a cuspidal automorphic representation.  This proves the first claim of Theorem \ref{iotads}.  
\item  Let $\tau$ be a representation of the $L$-group $\LGr$ and let $\CL(\Pi,\tau)$ be the corresponding local system on the curve $X$.  Then for any place $v$ of $K$, the local monodromy representation of the Galois group $\Gamma_v = Gal(K_v^{sep}/K_v)$ on $\CL(\Pi,\tau)$ is geometric, because $G$ is semisimple.
\item It follows by the properties of the local Genestier-Lafforgue correspondence -- specifically compatibility with parabolic induction -- that the exponents of the local representation $\Pi_v$ are Weil $q$-numbers.  In particular they are elements of an algebraic number field $E$ whose complex absolute values are independent of the embedding $\iota:  E \ra \bC$.
\item Suppose $\pi$ has a model $\pi_E$ over the number field $E$ and there is a complex embedding $\iota:  E \ra \bC$ such that 
the base change $\pi = \pi_E \otimes_{E,\iota} \bC$ is in the discrete series.   It then follows from Casselman's criterion that the exponents of $\pi$ are complex numbers satisfying certain inequalities.  By the previous observation, this remains true for any complex embedding of $E$.  This completes the proof of Theorem \ref{iotads}.  The proof of Theorem \ref{iotat} is similar.

\end{enumerate}

\begin{remark}  The analogue of Theorem \ref{iotads} is well known when $F$ is a $p$-adic field.  Clozel sketched an argument during the 1988 Ann Arbor conference, using the characterization of discrete series in terms of limit multiplicities in $\CA_0(G,U,\bC)$ when $G$ is anisotropic and the local component at a chosen place $v$ shrinks to the identity.  A complete proof of this argument was published by Dat in \cite[Prop. 4.7]{Dat}.  Probably a similar argument will be possible for $F$ of equal characteristic when there is a complete invariant trace formula over function fields.
\end{remark}


\subsection{The Aubert-Zelevinsky involution}

We let $F = k((t))$ be our local field.  We assume $F$ is an extension of our given $K$, so we can define the group $G(F)$ of $F$-rational points. 

\begin{prop}\label{AZcusp}  Let $\pi$ be a tempered (irreducible admissible) representation of $G(F)$, and let $AZ(\pi)$ denote the image of $\pi$ under the Aubert-Zelevinsky involution.  Then $\pi$ and $AZ(\pi)$ have the same cuspidal support, in the strong sense that $\pi$ and $AZ(\pi)$ are both constituents of the same (normalized) parabolically induced representation.

In particular, the semisimple local $\ell$-adic parameters of $\pi$ and $AZ(\pi)$ coincide:
$$\CL^{ss}(AZ(\pi)) = \CL^{ss}(\pi).$$
\end{prop}
\begin{proof}  The preservation of the cuspidal support under $AZ$ is immediate from \cite[Theorem 1.7]{Aubert95}, see also \cite[proof of Corollary 3.9]{Aubert95}.
\end{proof}

If $\pi$ is discrete series the semisimple parameter $\CL^{ss}(\pi)$ can be completed uniquely to a Weil-Deligne parameter $(\CL^{ss}(\pi),N_\pi)$ that satisfies purity of the monodromy weight filtration -- what is called Weil-Deligne purity in \cite{GHS}, where this is proved.   More generally, as we explain below, we can do the same for any tempered $\pi$.  In that case we expect that the Weil-Deligne parameter attached to $AZ(\pi)$ should just be the pair $(\CL^{ss}(\pi),0)$ -- in other words, the monodromy operator $N_{AZ(\pi)}$ is $0$ -- but we know of no global construction that can be used to confirm this expectation.

On the other hand, it is conjectured  that the {\it Arthur} $SL(2)$ attached to $AZ(\pi)$ should have a nilpotent operator $N_{Art}$ which is precisely the $N$ attached to $\pi$.  This is a special case of the expectation that the Aubert-Zelevinsky involution exchanges the Deligne $SL(2)$ of $\pi$ with the Arthur $SL(2)$ of $AZ(\pi)$, and vice versa.\footnote{This expectation was first stated as a precise conjecture in print in the articles \cite{Hi, Ban}, but it seems to have been in the air for some time before that.  For $G = GL(n)$ it had already been proved by Moeglin-Waldspurger and Tadi\'c, in articles cited by Ban.}

\section{The main theorems in the generic case}\label{mtgc}

Let $\ft^\vee$ denote a Cartan subalgebra of $\Lie(\hG)$, $\ft^\vee_\RR$ its real form.  Let $\sigma \in \Phi(U)$  be a parameter as in Theorem \ref{SeqT}, and
let $\pi \subset \CA_0(G,\bar{\QQ})$ be an irreducible representation that belongs to $\sigma$, as in Definition \ref{belongs}.  For any place $w \notin S(U)$ we let $s_w(\sigma)$ 
denote the component in $\hG$ of the local Satake parameter of $\pi$ at $w$.  We write 
\begin{equation}\label{logq}  \nu_w(\sigma) = \log_q(|\iota(s_w(\sigma)|) \in \ft^\vee_\RR.
\end{equation}

We record the following well-known fact.
\begin{lemma}\label{temperedprin}  Let $\sigma$ be as above.  Then the representation $\pi_w$ is tempered
if and only if $\nu_w(\sigma) = 0$.  Moreover, $\pi_w$ is then an irreducible summand of a principal series representation unitarily induced from a unitary character of the minimal
parabolic subgroup of $G(K_w)$.
\end{lemma}

\begin{proof}
By \cite[Corollary 7.2.2]{Cas}, we may assume without loss of generality that $\pi_w$ is an irreducible subrepresentation of a minimal spherical principal series. Then, if $\pi_w$ is a discrete series, the claim follows immediately from \cite[Theorem 6.5.1]{Cas} (with $\Theta=\emptyset$ in Casselman's notation). More generally, if $\pi_w$ is tempered, then $\pi_w$ can be embedded as a direct summand of a unitary generalized principal series induced from a discrete series on a Levi subgroup, see for example \cite[Theorem 3.22]{DO}, and the claim follows by induction in stages.
\end{proof}

For any irreducible $M$-dimensional representation $\tau$ of $\LGr$, the composition 
$$\tau\circ\sigma:  Gal(K^{sep}/K)~ \ra ~GL(M,\Qlb)$$
corresponds to a semisimple $\ell$-adic local system, denoted $\CL(\sigma,\tau)$, on the curve $X \setminus |S|$ over the finite field $k_1$.   We write
\begin{equation}\label{CLpt}
\CL(\pi,\tau) = \CL(\sigma,\tau)
\end{equation}
 if $\pi$ belongs to $\sigma$.    


We now state our main result in the generic case.

\begin{thm}\label{mainthm}  Let $\pi \subset \CA_0(G,\bar{\QQ})$ with $m_U(\pi) \neq 0$, and suppose
$\pi$ belongs to the parameter $\sigma$.  
Suppose that:
\begin{enumerate}
\item[(1)] the local component $\pi_v$ at our chosen {quasi-}split unramified place $v$ is generic.  
\item[(2)] there is a place $u \in S$ such that $\pi_u$ belongs to the discrete series.   
\end{enumerate}
Then $\pi_w$ is tempered for all $w \notin S$ (at which $G$ is unramified).

More generally, with $w$ and $u$ as above, suppose $\pi_u$ is tempered.  Then $\pi_w$ is tempered for all $w \notin S$.  In particular, if $v = u$
and $\pi_v$ is generic and tempered then  $\pi_w$ is tempered for all $w \notin S$.
\end{thm}

If $G$ is of type $C_n$ or $D_n$, we assume to begin with that  $G$ is adjoint, and  subsequently reduce to this case.   We may thus 
choose   representations $\tau_0$ of $\LGr$ so that the restriction to the Langlands dual group $G^\vee$
contains a representation with highest weight $\lambda^\vee$, determined by the following table:
 \begin{enumerate} 
\item[(i)] If $G=PGL(n)$, $n$ even, then $\lambda^\vee = \omega_1^\vee$, with the numbering of the roots in the Dynkin diagram:
\[\xymatrix{1\ar@{-}[r] &2 \ar@{-}[r] &3  \ar@{-}[r] &\dotsb \ar@{-}[r] &(n-1)};
\]
\item[(ii)] If $G=SO(2n+1)$, then $\lambda^\vee = \omega_n^\vee$ for the diagram
\[\xymatrix{1\ar@{<=}[r] &2 \ar@{-}[r] &3  \ar@{-}[r] &\dotsb \ar@{-}[r] &n};
\]
\item[(iii)] If $G=PSp(2n)$, then $\lambda^\vee = \omega_1^\vee$ for the diagram
\[\xymatrix{1\ar@{=>}[r] &2 \ar@{-}[r] &3  \ar@{-}[r] &\dotsb \ar@{-}[r] &n};
\]
\item[(iv)] If $G=PSO(2n)$, then we must take both $\lambda^\vee = \omega_1^\vee$ and $\lambda^\vee = \omega_n^\vee$ for the diagram
\[\xymatrix{1\ar@{-}[rd]\\ &3 \ar@{-}[r] \ar@{-}[ld] &4  \ar@{-}[r] &\dotsb \ar@{-}[r] &n;\\
2}
\]
\item[(v)] If $G=E_7$, then $\lambda^\vee = \omega_7^\vee$ for the diagram
\[\xymatrix{1\ar@{-}[r] &3 \ar@{-}[r]  &4\ar@{-}[d]  \ar@{-}[r] &5 \ar@{-}[r] &6 \ar@{-}[r] &7.\\
&&2}
\]
\end{enumerate}

 The relevant results regarding the unitarity of local factors are based on the classification of generic unitary spherical representations, which is recalled in \S \ref{classificationsection} and \S \ref{s:nonsplit}.
  
  \medskip

 \begin{proof}  First assume $\pi_u$ is in the discrete series.  Recall that we are assuming, temporarily, that $G$ is adjoint
 if it is of type $C_n$ or $D_n$.  Let $\tau$ be either the adjoint representation of $\LGr$ or the representation $\tau_0$ in the above list.
By \cite[Th\'eor\`eme (3.4.1) (i)]{De} or \cite[Corollary VII.8]{Laf02}, the semisimple $\ell$-adic local system $\CL(\pi,\tau)$ can be written as a (finite) direct sum
\begin{equation}\label{chidecomp} \CL(\pi,\tau) = \oplus_{\chi/\sim} \CL(\pi,\tau)_\chi \otimes \chi
\end{equation}
where each $\CL(\pi,\tau)_\chi$ is a punctually $\iota$-mixed local system with integer weights
and $\chi$ is the pullback to $C \setminus |S|$ of a rank $1$ $\ell$-adic local system over $Spec(k_1)$, in other words a continuous $\ell$-adic character of $Gal(\bar{k}_1/k_1)$.     Moreover, this decomposition is unique, when $\chi$ is taken to run over the equivalence classes for the following relation:  $\chi$ and $\chi'$ are equivalent, written $\chi \sim \chi'$ in \eqref{chidecomp}, if $\chi'\cdot \chi^{-1}$ is pure of some (integer) weight.

Consider the monodromy representation $\CL(\pi,\tau)_\chi\otimes \chi$ at $u$.    Let
$\CL_F(\pi_u)$ denote the semisimple $\ell$-adic local parameter attached to $\pi_u$  as in Theorem \ref{sslocal}.                   
Each irreducible composition factor of $\CL(\pi,\tau)_\chi\otimes \chi$ is then a  constituent of  $\tau\circ \CL_F(\pi_u)$.  
But it is proved in \cite{GHS} that $\tau\circ \CL_F(\pi_u)$
extends uniquely to a {pure Weil-Deligne representation} of weight $0$ for any $\tau$, in the sense of Definition \ref{mwf}.
By the uniqueness of the decomposition \eqref{chidecomp}, this implies that $\CL(\pi,\tau)$ is itself mixed.  

It follows that every irreducible component of  the semisimple $\ell$-adic local system $\CL(\pi,\tau)$ is punctually $\iota$-pure, for every $\iota$.   Thus for any
 $w \notin S$, the eigenvalues of $\tau\circ \sigma(\Frob_w)$ are Weil $q_w$-numbers, where $q_w$ is the order of the residue field at $w$.   This in turn implies that for any
 weight $\alpha$ of $\tau$, we have
 \begin{equation}\label{qnumbers} \langle \alpha,\nu_w(\sigma) \rangle \in \frac{1}{2}\ZZ.
 \end{equation}
 
In particular, this is the case when $w = v$, where $\pi_v$ is generic. 
Applying this when $\tau$ is the adjoint representation and the representation $\tau_0$ above, Corollary \ref{c:half-integral}(1) (when $G(K_v)$ is split) and Corollary \ref{c:unramified} (when $G(K_v)$ is quasi-split nonsplit)  then imply that $\nu_v(\sigma) = 0$. 
 It follows that, for every $\iota$, every irreducible component of $\CL(\pi,\tau)$ is punctually $\iota$-pure of weight $0$ at $v$, and therefore at every unramified place $w$.  Lemma \ref{temperedprin} then implies that $\pi_w$ is tempered at every unramified place $w$.

If $G$ is of type $C_n$ or $D_n$, we now drop the assumption that $G$ is adjoint.  At this point we would like to say that there is a finite (ramified) cover $X'$ of the curve $X$, with function
field $K'$ such that
$\pi$ admits a base change $\pi'$ to a cuspidal representation of   $G(\ad_{K'})$ with trivial central character.  Then $\pi'$ is a constituent of the pullback
to $G(\ad_{K'})$ of a cuspidal representation of the adjoint quotient of $G$, and we can argue as before.  Since the stable twisted trace formula is not currently available over
function fields, the existence of such a $\pi'$ is not known.  Instead we argue directly on the side of the parameters.  As above, we assume $\pi$ belongs to the $\ell$-adic parameter
$\sigma:  Gal(K^{sep}/K) \ra \LGr(\Qlb)$.  Now the obstruction to lifting $\sigma$ to a homomorphism $\tilde{\sigma}:  Gal(K^{sep}/K) \ra \LGad(\Qlb)$ lies in
$H^2(Gal(K^{sep}/K),Z)$, where $Z$ is the (finite) center of $\LGad(\Qlb)$.  It follows easily that there is a finite extension $K'$ of $K$, as above, such that
the restriction $\sigma'$ to $Gal(K^{sep}/K')$ lifts to a homomorphism $\tilde{\sigma'}:  Gal(K^{sep}/K') \ra \LGad(\Qlb)$.  Let $X'$ be the finite cover
of $X$ with function field $K'$, and let $S'$ be the set of places of $X'$ where $\tilde{\sigma'}$ is ramified.  Now we can compose with $\tau_0$ in the above
list to obtain a local system $\CL(\pi,\tau_0)$ over $X' \setminus |S'|$.  If $u'$ is a place of $K'$ above $u$, then the argument above shows that every component of the
restriction of $\CL(\pi,\tau_0)$ to the decomposition group at $u'$ extends to a pure Weil-Deligne representation of weight $0$, and we conclude as before.\footnote{By a well-known Baire category argument, we know that the image of $\sigma$ lies in $\LGr(E)$ for some finite extension $E$ of $\mathbb{Q}_\ell$.  Johan de Jong has pointed out that this implies that there is an open subgroup of the image of $\sigma$ that lifts to  $\LGad(\Qlb)$, so there is no need even to mention
the obstruction in Galois $H^2$.  }

Suppose now that $\pi_u$ is tempered.  Thus there is a parabolic subgroup $P \subset G(K_u)$ with Levi component $M$ and a discrete series representation $\rho$ of $M$, with unitary central character,
such that $\pi_u$ is an irreducible constituent of the  unitarily induced representation $\Ind_P^{G(F)} \rho$.  Since $\rho$ has unitary central character and
$\pi_u$ is a tempered local component of a cuspidal automorphic representation, it follows from Theorem \ref{iotat} that $\rho$ and $\pi_u$ are defined over a number field $E \subset \bC$, say $\rho = \rho_E \otimes_{E} \bC$; moreover, for any embedding $\iota:  E \ra \bC$, 
$$\rho_\iota = \rho_E\otimes_{E,\iota} \bC$$
has unitary central character.   Thus the Frobenius eigenvalues of  $\CL^{ss}(\rho)$ are Weil $q_u$-numbers,  
where $q_u$ is the order of the residue field at $u$.  
By the compatibility of the Genestier-Lafforgue correspondence with unitary parabolic induction, the same is true of the
Frobenius eigenvalues of $\CL^{ss}(\pi_u)$.  The remainder of the argument follows as in the discrete series case.  
\end{proof}

\begin{rmk*}  When $\pi_u$ is not supercuspidal, the proof in \cite{GHS} that $\tau\circ \CL^{ss}(\pi_u)$ is a pure Weil-Deligne representation of weight $0$ is based on an original argument of 
Beuzart-Plessis based in turn on the Deligne-Kazhdan simple trace formula, which is currently the most useful version of the trace formula generally available over function fields.  Beuzart-Plessis's argument is contained in the appendix to  \cite{GHS}.
\end{rmk*}

\begin{cor}\label{locgeneric}  Let $\pi$ be as in the statement of Theorem \ref{mainthm}.  Then $\pi_w$ is generic at all unramified places $w$  such that $G(K_w)$ is quasi-split.  In particular, if $G$ is a quasi-split group then $\pi$ is locally generic almost everywhere.
\end{cor}

\begin{proof}  It remains to record that if $\pi_w$ is spherical and tempered, then it is generic. This is well known, for example, for the general case, see \cite[Proposition 7.4]{Re}, the subtlety having to do with the isogeny class of $G(K_w)$. When $G(K_w)$ is split adjoint, the claim is easy, and for the benefit of the reader, we explain how it follows from \cite{KL}.  By the classification theorems of \cite{KL}, the smooth irreducible $G(K_w)$-representations with Iwahori fixed vectors are in one-to-one correspondence with  $G^\vee$-conjugacy classes of triples:
\[(s,e^\vee,\rho):\ s\in G^\vee \text{ semisimple},\ e^\vee\in \mathfrak g^\vee,\ \text{Ad}(s)e^\vee= q e^\vee,
\]
and $\rho$ is an irreducible representation of ``Springer type" of the group of components of the centralizer of $s$ and $e^\vee$ in $G^\vee$, such that the center of $G^\vee$ acts by the identity in $\rho$. Notice that $e^\vee$ must necessarily be nilpotent. Let $\pi_w(s,e^\vee,\rho)$ denote the corresponding $G(K_w)$-representation. Let $\{e^\vee,h^\vee,f^\vee\}$ be a Lie triple in $\mathfrak g^\vee$. We may arrange so that $s\in T^\vee$ (the maximal torus of $G^\vee$) and $h^\vee\in \mathfrak t^\vee_{\mathbb R}$. Set $s_0=s q^{-h^\vee/2}$ so that $\text{Ad}(s_0)e^\vee=e^\vee$. The Kazhdan-Lusztig correspondence has a number of properties including:
\begin{enumerate}
\item $\pi_w(s,e^\vee,\rho)$ is spherical if and only if $e^\vee=0$. Since the centralizer of $s$ in $G^\vee$ is connected, $G^\vee$ being simply connected, the representation $\rho$ is automatically trivial. Moreover, $s$ is the Satake parameter of $\pi_w$.
\item $\pi_w(s,e^\vee,\rho)$ is tempered if and only if $s_0$ is compact.
\end{enumerate}
In our case, $\pi_w$ is both spherical and tempered, therefore $\pi_w=\pi_w(s,0,\text{triv})$, where $s=s_0$ is compact. But then the $q$-eigenspace of $\text{Ad}(s)$ in $\mathfrak g^\vee$ is zero, so $\pi_w$ is the only irreducible representation with semisimple parameter $s$, hence it must be the full minimal spherical principal series with Satake parameter $s$. In particular, $\pi_w$ is generic.
\end{proof}

Shahidi's $L$-packet conjecture asserts that every tempered $L$-packet contains a globally generic member.
Since two members of a hypothetical $L$-packet are locally isomorphic almost everywhere, Corollary \ref{locgeneric} is consistent with Shahidi's conjecture.  

\begin{remark}  There are precise conjectures \cite[(8.3)]{LafICM} to the effect  that the following is always true: 

\begin{Expectation}\label{stack}  The irreducible components of the $\ell$-adic local systems $\CL(\pi,\tau) = \CL(\sigma,\tau)$ are realized in the total direct image $Rf_!(\Qlb)$ with compact support, where
$f:  Z \ra C \setminus |S'|$ is a morphism of algebraic stacks 
of finite type
for some finite set $S' \supset S$, when $\tau$ is the adjoint representation of $\LGr$ on $\Lie(\hG)$ and also when $\tau = \tau_0$ (if $G$ is not both simply connected and adjoint).
\end{Expectation}

The stacks in question are precisely the Deligne-Mumford
moduli stacks of shtukas that figure in Lafforgue's constructions of Galois parameters.
Assuming this to be the case, it is well-known to experts that every irreducible component of  the semisimple $\ell$-adic local system $\CL(\pi,\tau)$ is punctually $\iota$-pure with integer weights.  However, there seems to be no accessible proof in the literature.  We therefore provide a sketch of a proof in  \S \ref{purity}.  
Admitting the purity for the moment we see that, for any
 $w \notin S(U)$, the eigenvalues of $\tau\circ \sigma(\Frob_w)$ are Weil $q_w$-numbers, where $q_w$ is the order of the residue field at $w$.   This implies as above that  any
 weight $\alpha$ of $\tau$ satisfies \eqref{qnumbers}.  In particular, this is the case when $w = v$, where $\pi_v$ is generic.
Reasoning as in the proof of Theorem \ref{mainthm}, Corollary \ref{c:half-integral}  then implies that $\nu_v(\sigma) = 0$. 

 It follows that every irreducible component of $\CL(\pi,\tau)$ is punctually $\iota$-pure of weight $0$ at $v$, and therefore at every unramified place $w$.  Lemma \ref{temperedprin} then implies that $\pi_w$ is tempered at every unramified place $w$.
 
 Lafforgue and Zhu have proved Expectation \ref{stack} to be true when the global parameter of $\pi$ is {\it elliptic}
 \cite[Proposition 1.2]{LZ}.   Presumably any $\pi$ satisfying the hypotheses of Theorem \ref{mainthm} is elliptic.  For elliptic parameters, however, Theorem \ref{mainthm} is superfluous; it is explained in \S 3 of \cite{LZ} that such $\pi$ are already tempered at  unramified places.   We thank Zhu for spelling out the argument, which goes roughly
 as follows:  letting $M$ denote the Zariski closure of the image in ${}^LG$ of the (semisimple) parameter of $\pi$, if the parameter is not tempered then $M$ has non-trivial characters, corresponding to the weight decomposition of the semisimplification of $\CL(\pi,\tau)$.  But this implies the centralizer of the parameter of $\pi$ is of positive dimension, hence the parameter  is not elliptic by definition.
 \end{remark}
  

\begin{sketch}[purity]\label{purity} We thank Martin Olsson for the following sketch of the proof of the purity claim.  It suffices in fact to show that the cohomology $Rf_*(\Qlb)$ is $\iota$-mixed.   The main reference is Theorem 2.11
 of the paper \cite{Su}, which proves that the property of being $\iota$-mixed is preserved under the six operations of Laszlo and Olsson (the stack version of Grothendieck's six operations) and Verdier duality.  We thus need to show that $Rf_!$ preserves the property that the weights are integral.   Now one can show that Verdier duality preserves the integrality by using a smooth cover to reduce to the case of schemes, so we may replace $Rf_!$ by $Rf_*$.  As in the argument in \cite{Su}, one reduces to the case when the target stack is a scheme. 
  Thus it suffices to consider the case $f: X \ra Y$, with $X$ an algebraic stack and $Y$ a scheme.  By using a hypercover of $X$ by schemes and the associated spectral sequence we are thus reduced to the case of schemes.  By Verdier duality again, this reduces to the case of $Rf_!$, which is proved in \cite{De}.
\end{sketch}

\begin{remark}  The proof of Theorem \ref{mainthm} implies, by the arguments in \cite{GHS}, that the semisimple
Genestier-Lafforgue parameters $\CL^{\rm ss}(\pi_u)$ of the local components $\pi_u$ of $\pi$ at all places, including $u \in S$, extend (uniquely) to tempered Weil-Deligne parameters, in the sense defined in \cite{GHS}.  But this does not imply that the ramified components $\pi_u$ themselves are necessarily tempered, although it is generally believed that they must be.

\end{remark}

\section{Half-integral generic spherical unitary points}\label{classificationsection}

We fix notation. Let $k$ be a nonarchimedean local field with discrete valuation $\val_k$, ring of integers $\mathfrak O$, and finite residue field $\bF_q$. Let $G$ be a quasi-simple split $k$-group of adjoint type and $K=G(\mathfrak O)$, the maximal compact hyperspecial subgroup of $G(k)$. Let $B\supset T$ denote a $k$-rational Borel subgroup and a maximal $k$-split torus, respectively. Denote the corresponding based root datum of $(G,B,T)$ by $(X,\Phi, X^\vee,\Phi^\vee,\Pi)$, where $\Pi$ are the simple roots, and let $W$ be the finite Weyl group. Let $\Phi^+\supset\Pi$ denote the positive roots, and $\Phi^{\vee,+}$ the corresponding positive coroots. Let $\langle ~,~\rangle$ be the natural pairing between $X$ and $X^\vee$. Let $G^\vee$ be the complex Langlands dual group with maximal torus $T^\vee=X\otimes_\bZ \bC^\times$. Recall the polar decomposition $T^\vee=T^\vee_\bR\cdot T^\vee_c$, where $T^\vee_\bR=X\otimes_\bZ \bR_{>0}$, $T^\vee_c=X\otimes_\bZ S^1$. 

If $s\in T^\vee\simeq\Hom(X^\vee,\bC^\times)$, let $\chi_s\in X$ denote the unramified character of $T(k)$,
\begin{equation}
\chi_s=s\circ \val_T: T(k)\to \bC^\times,
\end{equation} 
where  $\val_T:T(k)\to X^\vee$ is 
defined by 
\[\langle \val_T(t),\lambda\rangle=\val_k(\lambda(t)),\quad \text{for all }\lambda\in X\text{ and all }t\in T(k).
\]
Regard $\chi_s$ as a character of $B(k)$ by pullback from $T(k)$ and let
\[X(s)=i_{B(k)}^{G(k)}(\chi_s)
\]
be the (unitarily) induced spherical minimal principal series. The representation $X(s)$ has a unique irreducible $K$-spherical subquotient $L(s)$, and there is a one-to-one correspondence between irreducible $K$-spherical $G(k)$-representations and $W$-orbits in $T^\vee$:
\[s\leftrightarrow L(s).
\]
The semisimple element $s$ (or rather its $W$-orbit in $T^\vee$) is called the \emph{Satake parameter} of the irreducible $K$-spherical representation $L(s)$. 

Let $\nu\in \ft^\vee_\bR=\Lie(T^\vee_\bR)$ and write
\[s_\nu=\exp(\nu\log(q))\in T^\vee_\bR.
\]
Following \cite{BM}, one may refer to these parameters as \emph{real Satake parameters}. Every Satake parameter $s\in T^\vee$ can be written uniquely as
\begin{equation}\label{e:polar}
s=s_c\cdot s_\nu,\text{ for some } s_c\in T^\vee_c\text{ and }\nu\in  \ft^\vee_\bR.
\end{equation}

\begin{defn}[cf. Definition \ref{algebraic}] Let $(V,\rho)$ be a finite-dimensional $\mathfrak g^\vee$-representation. An element $\nu\in \ft^\vee_\bR$  is called \emph{half-integral} with respect to $V$, if for all weights $\chi$ of $V$, $\langle\chi,\nu\rangle\in \frac 12\mathbb Z$. In this case, we call the Satake parameter $s=s_c\cdot s_\nu$ \emph{geometric} (with respect to $V$).
\end{defn}

For simplicity, we also write
\[X(\nu):=X(s_\nu),\quad \nu \in \ft^\vee_\bR,
\]
and $L(\nu)$ for its irreducible $K$-spherical subquotient. (When $\nu$ is dominant, $L(\nu)$ is the unique irreducible quotient of $X(\nu)$.) With this notation, the trivial $G(k)$-representation is $L(\rho)$, where $\rho=\frac 12\sum_{\alpha\in \Phi^+}\alpha\in \ft^\vee_\bR$.

The irreducible unitary spherical $G(k)$-representations with real Satake parameters are classified in \cite{Ba} for symplectic and orthogonal split groups, \cite{Ci-F4,Ci-E6,BC-E8} for exceptional split groups, via solving the equivalent corresponding problem for Iwahori-Hecke algebras \cite{BM}. The case $G=GL(n)$ is by now well known and it is due to Tadi\' c \cite{Ta}.

By \cite{BM-generic}, the $K$-spherical quotient $L(s)$ is generic (in the sense of admitting Whittaker vectors) if and only if $X(s)=L(s)$, that is, if and only if $X(s)$ is irreducible. The reducibility points of $X(s)$ are well known, going back to Casselman \cite[Proposition 3.5]{Cas-spherical} for $s$ regular. When $\nu\in  \ft^\vee_\bR$, we have (see  \cite[\S3.1]{BC-E8}), 
\[X(\nu)=L(\nu) \text{ if and only if } \langle\al^\vee,\nu\rangle\neq 1\text{ for all }\al^\vee\in \Phi^{\vee,+}.
\]

Denote
\begin{equation}
\begin{aligned}
\CC_0&=\{\nu \in  \ft^\vee_\bR\mid \langle \al^\vee,\nu\rangle\ge 0,\text{ for all }\al\in \Pi\},\\
\quad \CU^g_0&=\{\nu\in \CC_0\mid   \langle\al^\vee,\nu\rangle\neq 1\text{ for all }\al^\vee\in \Phi^{\vee,+}, \text{ and }L(\nu)\text{ is unitary}\}.
\end{aligned}
\end{equation}
By the previous discussion $\CU^g_0$ parameterizes the (equivalence classes of) spherical generic representations with real Satake parameter. Moreover, for an irreducible spherical representation $L(\nu)$ to admit a nonzero invariant Hermitian form it is necessary and sufficient (\cite[\S3.1-3.2]{BC-E8}) that 
\begin{equation}
w_0(\nu)=-\nu,\text{ where } w_0 \text{ is the longest Weyl group element}.
\end{equation}
This is automatic when $\Phi^\vee$ is of types $B,C, G_2, F_4, E_7, E_8$. Denote 
\[\CC_{0,h}=\{\nu\in \CC^\vee_0\mid w_0(\nu)=-\nu\},
\]
so that $\CU^g_0\subset \CC_{0,h}.$

The explicit description of $\CU^g_0$ is given in {\it loc. cit.}. The nature of the answer is the following. The arrangement of hyperplanes
\[\alpha^\vee=1,\quad \alpha^\vee\in \Phi^{\vee,+},
\]
partitions the  fundamental Weyl chamber $\CC_0$ into $m$ open regions, 
\[m=\prod_{i=1}^\ell \frac{d_i+h}{d_i},\]
where $h$ is the Coxeter number, $\ell$ is the rank, and $d_i$ are the fundamental degrees for the \emph{coroot system} $\Phi^\vee$. For example, for $F_4$ there are $105$, while for $E_8$, there are $25080$ such open regions. It is clear that $\CU^g_0$ is a union of open regions in this arrangement of hyperplanes intersected with $\CC_{0,h}$. 

For each simple root system, denote by $\al_i$, $\al_i^\vee$, $\omega_i$, $i=1,\dots,\ell$ the simple roots, simple coroots, and fundamental weights, respectively, and
\[\rho=\sum_{i=1}^\ell \om_i.
\]
Consider the poset of positive coroots. We say that a positive coroot $\al^\vee$ has level $r$ if $\langle \al^\vee,\rho\rangle=r$. Let $\gamma^\vee$ be the highest positive coroot, the unique coroot of level $h-1$. Recall that an \emph{alcove} is a connected (open) component in 
\[
\ft^\vee_\bR\setminus \bigcup_{\alpha^\vee\in\Phi^\vee, m\in \bZ_{\ge 0}} \{\nu\in \ft^\vee_\bR \mid \langle\alpha^\vee,\nu\rangle=m\}.
\]
In particular, the \emph{fundamental alcove} is
\begin{equation}
\CA_0=\{\nu\in \CC_0\mid \langle \gamma^\vee,\nu\rangle<1\}.
\end{equation}

\begin{lemma}
$\CA_0\cap \CC_{0,h}\subseteq \CU^g_0.$
\end{lemma} 

\begin{proof}
At $\nu=0$, the spherical principal series $X(0)$ is unitary and irreducible. By a classical deformation argument (\emph{complementary series}) originally used in \cite[Theorem 3.3]{KS}, if we deform $\nu\in \CC_{0,h}$ continuously, $X(\nu)$ is unitary as long as it remains irreducible. For $\nu\in \CC_0$, the closest reducibility hyperplane to $0$ is $\gamma^\vee=1$.  This implies that $\CA_0\cap  \CC_{0,h}$ is unitary.
\end{proof}

Another classical result says that no unbounded region can be unitary, see \cite{HM} (or \cite[\S3.3]{BC-gen} for a proof in this setting):

\begin{lemma}\label{l:rho}
If $\nu\in \CC_0$ is such that $\langle\al^\vee_i,\nu\rangle>1$ for a simple coroot $\al^\vee_i$, then $\nu\notin \CU^g_0.$
\end{lemma}

But beyond these classical results, the determination of $\CU^g_0$ is difficult and has been achieved via a case-by-case analysis. We will use some elements of the explicit description which will be recorded in the sequel, but here let us give the general form of the answer.

\begin{thm}[\cite{Ba,Ci-F4,Ci-E6,BC-E8}]\label{t:alcoves} Suppose $G$ is a quasi-simple split group of adjoint type. The set $\CU_0^g$ is the intersection of the hermitian locus $\CC_0$ with a union of $2^d$ alcoves in $\CC_0$ (including the fundamental alcove $\CA_0$), where
\begin{itemize}
\item  $d=0$, if $G=PGL(n)$ or $SO(2n+1)$;
\item $d=\lfloor\frac {n-1}2\rfloor$ is $G=PSp(2n)$;
\item $d=n-1$ if $G=PSO(4n)$, $n\ge 2$, or if $G=PSO(4n+2)$, $n\ge 1$;
\item $d=1$, if $G=G_2,$ $F_4$, or $E_6$;
\item $d=3$, if $G=E_7$;
\item $d=4$, if $G=E_8$.
\end{itemize}
\end{thm}

\subsection{Integral unitary points for the adjoint $G^\vee$-representation}

We prove first that there are no nonzero generic unitary parameters $\nu$ that are integral with respect to the adjoint $G^\vee$-representation. 

\begin{prop}\label{p:integral-unitary}
If $\nu\in\CC_0$ is such that $\langle\alpha^\vee,\nu\rangle\in \bZ_{>0}$ for a positive coroot $\alpha^\vee$, then $\nu\notin \CU^g_0.$
\end{prop}

\begin{proof}
Let $\nu\in\CC_{0,h}$, otherwise it can't be a unitary parameter since it is not hermitian. If $\langle\alpha^\vee,\nu\rangle=m$ for a positive coroot $\alpha^\vee$ and a positive integer $m$, then by definition, $\nu$ does not belong to any alcove in $\CC_0$. Since $\CU_0^g$ is a union of alcoves by Theorem \ref{t:alcoves}, $\nu\notin \CU_0^g$.
\end{proof}

\subsection{Half-integral unitary points for the adjoint $G^\vee$-representation}

We would like to determine the generic unitary parameters $\nu$ that are also $\epsilon$-half-integral with respect to the adjoint $G^\vee$-representation for a fixed real number $0\le\epsilon<\frac 12$. In other words, we need to calculate:
\begin{equation}
\CU^{g,\frac 12+\epsilon}_0=\CU^g_0\cap \CC_0^{\frac 12+\epsilon},\text{ where }\CC_0^{\frac 12+\epsilon}=\{\nu\in \CC_0\mid \langle \al^\vee,\nu\rangle\in \epsilon+\frac 12\bZ,\text{ for all }\al^\vee\in \Phi^{\vee}\}.
\end{equation}
Since the condition has to be satisfied for both $\al^\vee$ and $-\al^\vee$, it follows immediately that $2\epsilon\in \frac 12\bZ$, hence $\CC_0^{\frac 12+\epsilon}=\emptyset$ unless 
\[\epsilon\in \{0,\frac 14\}.
\]
If the $\Phi$ contains a simple root subsystem of rank at least $2$ (in particular, if $G$ is quasi-simple and $G\neq PGL(2)$), then there exist simple coroots $\al_1^\vee$ and $\al_2^\vee$ such that $\beta^\vee=\al_1^\vee+\al_2^\vee$ is also a coroot. Then 
\[\langle\beta^\vee,\nu\rangle=\langle \al_1^\vee,\nu\rangle+\langle \al_2^\vee,\nu\rangle\in 2\epsilon+\frac 12\bZ.
\]
On the other hand, $\langle\beta^\vee,\nu\rangle\in \epsilon+\frac 12\bZ$. From this, $\epsilon\in\frac 12\bZ$ and so $\epsilon=0$ necessarily. Hence, if $G$ is quasi-simple of rank at least $2$,
\begin{equation}
\CC^{\frac 12+\epsilon}_0=\emptyset\text{ unless }\epsilon=0.
\end{equation}
For $G=PGL(2)$, we will only be interested in $\epsilon=0$ as well. 

\medskip

We now consider the case $\ep=0$. The ordering of the roots and the explicit coordinates are given in the subsections, for exceptional groups they are the Bourbaki ones; the coordinates are such that the pairings are just the dot product of vectors. Write a point $\nu\in \CC_0$ in the \emph{fundamental weight coordinates}: 
\[\nu=\sum_{i=1}^\ell \nu_i\omega_i,\quad \nu_i\ge 0,\text{ for all }1\le i\le \ell.\]
Then $\nu\in\CC_0^{\frac 12}$ if and only if $\nu_i\in \frac 12\bZ_{\ge 0}$ for all $i$. By Lemma \ref{l:rho}, it follows that a necessary condition for $\nu\in \CU^{g,\frac 12}_0$ is
\begin{equation}\label{e:a-i}
\nu_i\in \{0,\frac 12\},\text{ for all } 1\le i\le \ell.
\end{equation}

\begin{thm}\label{t:half-adjoint}
Suppose $G$ is a quasi-simple $k$-split group of adjoint type. The set $\CU^{g,\frac 12}_0$ is as follows:
\begin{itemize}
\item $G=PGL(n)$, $n$ even, $\CU^{g,\frac 12}_0=\{0,\frac 12\omega_{\frac n2}\}$;
\item $G=SO(2n+1)$,  $\CU^{g,\frac 12}_0=\{0,\frac 12\omega_1\}$;
\item $G=PSp(2n)$,  $\CU^{g,\frac 12}_0=\{0,\frac 12\omega_n\}$;
\item $G=PSO(2n)$, $n$ even,  $\CU^{g,\frac 12}_0=\{0,\frac 12\omega_1,\frac 12\omega_2,\frac 12\omega_n\}$;
\item $G=PSO(2n)$, $n$ odd,  $\CU^{g,\frac 12}_0=\{0,\frac 12\omega_n\}$;
 \item $G=E_7$, $\CU^{g,\frac 12}_0=\{0,\frac 12\omega_7\}$;
 \item For all other cases, $\CU^{g,\frac 12}_0=\{0\}$.
 \end{itemize}
\end{thm}

\begin{proof}
The strategy for each simple root system is the following:

\begin{enumerate}
\item The classification of the generic spherical unitary dual implies, in particular, that there exists a maximum level $r_0\ge 1$ such that $\nu\in \CU^g_0$ \emph{only if} $\langle \al^\vee,\nu\rangle<1$ for \emph{all} $\al^\vee$ of level $r_0$. The fact that $r_0\ge 1$ is equivalent with Lemma \ref{l:rho}.

Record which points $\nu\in\CC_0^{\frac 12}$ satisfy $\langle \al^\vee,\nu\rangle<1$ for  \emph{all} $\al^\vee$ of level $r_0$. Call this set $\CC_0^{\frac 12,\le r_0}$.

\item For each one of the resulting points $\nu\in \CC_0^{\frac 12,\le r_0}$, check if $\nu$ belongs to one of the hyperplanes of reducibility $\alpha^\vee=1$, for a positive coroot $\alpha^\vee$. In all cases, it turns out that $0\neq \nu\in \CC_0^{\frac 12,\le r_0}$ either  lies on the reducibility hyperplane $\gamma^\vee=1$, where $\gamma^\vee$ is the highest coroot, or $\nu$ is one of the exceptions listed in the theorem. 

Finally, for each one of the listed exceptions, verify that the parameter $\nu$ lies in the one of the unitary regions of $\CU^g_0$; it turns out they are always in the fundamental alcove $\gamma^\vee<1$.
\end{enumerate}
The details for each simple root system are presented in the next subsections. 
\end{proof}

\subsection{$PGL(n)$} For type $A$, it is more convenient to use standard coordinates, rather than coordinates with respect to the fundamental weights. The positive roots are $\ep_i-\ep_j$, $1\le j\le i\le n$, and we write a typical dominant hermitian parameter $\nu\in \CC_0$ as:
\begin{align*}
\nu&=(a_1,a_2,\dots,a_k,-a_k,-a_{k-1},\dots,-a_1), &\text{ if }n&=2k,\\
\nu&=(a_1,a_2,\dots,a_k,0,-a_k,-a_{k-1},\dots,-a_1), &\text{ if }n&=2k+1,
\end{align*}
where $a_1\ge a_2\ge\dots\ge a_k\ge 0.$ The unitarity condition is \cite{Ta}:
\begin{equation}\label{bound-A}
0\le a_k\le a_{k-1}\le \dots\le a_2\le a_1< \frac 12.
\end{equation}
Suppose $n=2k$ (even). Then $\nu$ is half-integral if and only if $a_k\in \frac 14\bZ$ and $\nu_i-\nu_{i+1}\in \frac 12\bZ$, $1\le i\le k-1$. Hence the only possible nonzero point in the unitary region is
\[\frac 12\om_{k}, \text{ for which }\nu_1=\nu_2=\dots=\nu_k=\frac 14.
\]
This point is in the unitary region obviously.

\smallskip

If $n=2k+1$ (odd), then $\nu$ is half-integral if and only if $a_i\in \frac 12\bZ$, but no such point is in the unitary region (\ref{bound-A}).

\subsection{$SO(2n+1)$} In the coordinates $\ep_1,\dots,\ep_n$ in $\bR^n$,  we use the coroots 
\[\al_1^\vee=2\ep_1,\ \al_2^\vee=-\ep_1+\ep_2,\  \al_3^\vee=-\ep_2+\ep_3,\dots,\  \al_n^\vee=-\ep_{n-1}+\ep_n,
\]
and weights
\[\om_1=\frac 12(\ep_1+\ep_2+\dots+\ep_n),\quad \om_i=\ep_i+\ep_{i+1}+\dots+\ep_n,\ 2\le i\le n.
\]
The highest coroot is $\gamma^\vee=2\ep_n$. From \cite[Theorem 3.1]{Ba}, $r_0=2n-1$ and $\gamma^\vee$ is the only coroot of that level. If we write $\nu=\sum_{i=1}^n\nu_i\om_i \in\CC_0^{\frac 12}$, $\nu_1,\nu_2\in \frac 12\bZ_{\ge 0}$, then the unitary bound implies
\[\nu_1+2\nu_2+2\nu_3+\dots+2\nu_n<1,
\]
hence $\nu_2=\nu_3=\dots=\nu_n=0$. We have two points left: $0$ and $\frac 12\om_1$. Since $\langle \gamma^\vee,\frac 12\om_1\rangle=\frac 12<1$, this point is in the fundamental alcove, hence unitary.

\subsection{$PSp(2n)$}\label{sub:psp}  In the coordinates $\ep_1,\dots,\ep_n$ in $\bR^n$,  we use the coroots 
\[\al_1^\vee=\ep_1,\ \al_2^\vee=-\ep_1+\ep_2,\  \al_3^\vee=-\ep_2+\ep_3,\dots,\  \al_n^\vee=-\ep_{n-1}+\ep_n,
\]
and weights
\[\om_i=\ep_i+\ep_{i+1}+\dots+\ep_n,\ 1\le i\le n.
\]
The highest coroot is $\gamma^\vee=\ep_{n-1}+\ep_n$. From \cite[Theorem 3.1]{Ba}, we can deduce that 
\[r_0=\begin{cases} n,& n \text{ odd},\\n+1, &n \text{ even}.\end{cases}
\]
If $n$ is odd, the coroots at level $r_0$ are
\begin{equation}
\ep_1+\ep_{n-1},\ \ep_2+\ep_{n-2},\dots,\ \ep_{\frac {n-1}2}+\ep_{\frac {n+1}2},\ \ep_n,
\end{equation}
while if $n$ is even, they are
\begin{equation}
\ep_1+\ep_{n},\ \ep_2+\ep_{n-1},\dots,\ \ep_{\frac {n}2}+\ep_{\frac {n}2+1}.
\end{equation}
Write $\nu=\sum_{i=1}^n\nu_i\om_i \CC_0^{\frac 12}$, $\nu_1,\nu_2\in \frac 12\bZ_{\ge 0}$. Then these unitary bounds imply, when $n$ is odd:
\begin{align*}
2\nu_1+\nu_2+\nu_3+\dots+\nu_{n-1}&<1,\\
2\nu_1+2\nu_2+\nu_3+\dots+\nu_{n-2}&<1,\\
\vdots\\
2\nu_1+2\nu_2+\dots+2\nu_{\frac{n-1}2}+\nu_{\frac{n+1}2}&<1,\\
\nu_1+\nu_2+\dots+\nu_n&<1,
\end{align*}
and when $n$ is even:
\begin{align*}
2\nu_1+\nu_2+\nu_3+\dots+\nu_{n}&<1,\\
2\nu_1+2\nu_2+\nu_3+\dots+\nu_{n-1}&<1,\\
\vdots\\
2\nu_1+2\nu_2+\dots+2\nu_{\frac{n}2}+\nu_{\frac{n}2+1}&<1.\\
\end{align*}
In both cases, it follows immediately that
\[\nu_1=\nu_2=\dots=\nu_{\lfloor\frac n2\rfloor}=0,
\]
and
\[\nu_{\lfloor\frac n2\rfloor+1}+\nu_{\lfloor\frac n2\rfloor+2}+\dots+\nu_n<1.
\]
This means that at most one $\nu_j=\frac 12$, $j=\lfloor\frac n2\rfloor+1,\dots,n$, and the others are $0$. This gives the half-integral points $0$ and
\[\frac 12\om_j,\quad j=\lfloor\frac n2\rfloor+1,\dots,n.
\]
But for $j<n$, $\frac 12\om_j$ lies on the hyperplane $\gamma^\vee=1$, so the only nonzero point left is $\frac 12\om_n$. Since $\langle \gamma^\vee,\frac 12\om_n\rangle=\frac 12<1$, this point is in the fundamental alcove, hence unitary.

\subsection{$PSO(2n)$} In the coordinates $\ep_1,\dots,\ep_n$ in $\bR^n$,  we use the coroots 
\[\al_1^\vee=\ep_1+\ep_2,\ \al_2^\vee=-\ep_1+\ep_2,\  \al_3^\vee=-\ep_2+\ep_3,\dots,\  \al_n^\vee=-\ep_{n-1}+\ep_n,
\]
and weights
\[\om_1=\frac 12(\ep_1+\ep_2+\dots+\ep_n),\ \om_2=\frac 12(-\ep_1+\ep_2+\dots+\ep_n),\quad \om_i=\ep_i+\ep_{i+1}+\dots+\ep_n,\ 3\le i\le n.\]
The highest coroot is $\gamma^\vee=\ep_{n-1}+\ep_n$. Write $\nu=\sum_{i=1}^n\nu_i\om_i \CC_0^{\frac 12}$, $\nu_1,\nu_2\in \frac 12\bZ_{\ge 0}$. If $n$ is odd, we need to impose in addition that $\nu_1=\nu_2$ in order for this parameter to correspond to a hermitian spherical representation.

Suppose $n$ is even. From \cite[Theorem 3.1]{Ba}, we can deduce that $r_0=n-1$ and the coroots at this level are
\begin{equation}
\ep_1+\ep_n,\ -\ep_1+\ep_n,\ \ep_2+\ep_{n-1},\ \ep_3+\ep_{n-2},\dots,\ \ep_{\frac n2}+\ep_{\frac n2+1}.
\end{equation}
The unitary bound given by these coroots implies (via a similar discussion to that for $Sp(2n)$, $n$ even), that
\[\nu_3=\nu_4=\dots=\nu_{\frac n2}=0,
\]
and at most one of $\nu_1,\nu_2,\nu_j$, $j=\frac n2+1,\dots,n$, can equal $\frac 12$, the rest being zero. In addition to the origin $0$, we are left with
\[\frac 12\om_1,\frac 12\om_2,\quad \frac 12\om_j,\ j=\frac n2+1,\dots,n.
\]
All $\frac 12\om_j$, $\frac n2+1\le j\le n-1$ lie on the hyperplane $\gamma^\vee=1$, so the only nonzero points remaining are $\frac 12\om_1,\frac 12\om_2,\frac 12\om_n$, for which we have
\[\langle \gamma^\vee,\frac 12\om_1\rangle=\langle \gamma^\vee,\frac 12\om_2\rangle=\langle \gamma^\vee,\frac 12\om_n\rangle=\frac 12<1,
\]
meaning they are in the fundamental alcove and hence unitary.

\smallskip

If $n$ is odd, a similar analysis holds, except that now $\frac 12\om_1$ and $\frac 12\om_2$ are not hermitian.

\subsection{$G_2$} In Bourbaki coordinates, realize the root system of type $G_2$ in the hyperplane perpendicular to the vector $(1,1,1)$ in $\bR^3$:
\begin{equation}
\begin{aligned}
\al_1&=(2,-1,-1) & \al_1^\vee&=(\frac 23,-\frac 13,-\frac 13) &\omega_1&=(1,1,-2)\\
\al_2&=(-1,1,0) & \al_2 ^\vee&=(-1,1,0) &\omega_2&=(0,1,-1).
\end{aligned}
\end{equation}
Write $\nu=\nu_1\omega_1+\nu_2\omega_2\in \CC_0^{\frac 12}$, $\nu_1,\nu_2\in \frac 12\bZ_{\ge 0}$. From \cite[Appendix B]{Ci-F4}, we see that $r_0=3$ and there is only one coroot $\beta^\vee=2\al_1^\vee+\alpha_2^\vee$ at this level. Hence a necessary condition for unitarity is
\begin{equation}
2\nu_1+\nu_2<1.
\end{equation}
From this, $2\nu_1<1$, so necessarily $\nu_1=0$. Then $\nu_2<1$, so $\nu_2=0$ or $\nu_2=\frac 12$. If $\nu_2=\frac 12$, we get the point $\nu=\frac 12\omega_2$, but this lies on the hyperplane $\gamma^\vee=1$, where $\gamma^\vee=3\alpha_1^\vee+2\alpha_2^\vee$ is the longest coroot.

\subsection{$F_4$} In coordinates in $\bR^4$:
\begin{equation}
\begin{aligned}
\al_1&=\frac 12(1,-1,-1,-1) &\al_1^\vee&=(1,-1,-1,-1) &\omega_1&=(1,0,0,0)\\
\al_2&=(0,0,0,1) &\al_2^\vee&=(0,0,0,2) &\omega_2&=(\frac 32,\frac 12,\frac 12,\frac 12)\\
\al_3&=(0,0,1,-1) &\al_3^\vee&=(0,0,1,-1) &\omega_3&=(2,1,1,0)\\
\al_4&=(0,1,-1,0) &\al_4^\vee&=(0,1,-1,0) &\omega_4&=(1,1,0,0).
\end{aligned}
\end{equation}
From \cite[Table, p. 126]{Ci-F4}, we see that $r_0=9$ and there is only one coroot $\beta^\vee=(1,1,1,-1)=\al_1^\vee+2\al_2^\vee+4\al_3^\vee+2\al_4^\vee$ of level $9$.

Write $\nu=\sum_{i=1}^4\nu_i\omega_i\in \CC_0^{\frac 12}$, $\nu_i\in \frac 12\bZ_{\ge 0}$ for all $i$. Using $\langle\beta^\vee,\nu\rangle<1$, we get
\begin{equation}
\nu_1+2\nu_2+4\nu_3+2\nu_4<1.
\end{equation}
Immediately, $\nu_2=\nu_3=\nu_4=0$ and $\nu_1<1$. So either $\nu_1=\frac 12$ or $\nu_1=0$. In the first case, we get the point $\nu=\frac 12\omega_1$, but this lies on the hyperplane $\gamma^\vee=1$, where $\gamma^\vee=(2,0,0,0)$ is the longest coroot. Thus we are left with the origin.

\subsection{$E_8$} In coordinates in $\bR^8$ with respect to a basis $\ep_i$, $1\le i\le 8$:
\begin{equation}
\begin{aligned}
\al_1^\vee&=\frac 12(1,-1,-1,-1,-1,-1,-1,1), &\al_2^\vee&=\ep_1+\ep_2,& \al_3^\vee&=-\ep_1+\ep_2, & \al_4^\vee&=-\ep_2+\ep_3, \\
 \al_5^\vee&=-\ep_3+\ep_4,  &\al_6^\vee&=-\ep_4+\ep_5,&  \al_7^\vee&=-\ep_5+\ep_6,  &\al_8^\vee&=-\ep_6+\ep_7;\\
\end{aligned}
\end{equation}
\begin{equation}
\begin{aligned}
\om_1&=2\ep_8, &\om_2&=(\frac 12,\frac 12,\frac 12,\frac 12,\frac 12,\frac 12,\frac 12,\frac 52), &\om_3&=(-\frac 12,\frac 12,\frac 12,\frac 12,\frac 12,\frac 12,\frac 12,\frac 72),\\
 \om_4&=(0,0,1,1,1,1,1,5),& \om_5&=(0,0,0,1,1,1,1,4), &\om_6&=(0,0,0,0,1,1,1,3),\\
 \om_7&=(0,0,0,0,0,1,1,2), &\om_8&=\ep_7+\ep_8.
\end{aligned}
\end{equation}
From \cite[\S7.2.3]{BC-E8}, we see that $r_0=15$ and there are four coroots at level $15$. They are:
\begin{equation}\label{bound-E8}
\begin{aligned}
\beta_1^\vee&=\al_1^\vee+2\al_2^\vee+2\al_3^\vee+4\al_4^\vee+3\al_5^\vee+2\al_6^\vee+\al_7^\vee=\frac 12(1,-1,1,1,1,1,-1,1) \\
\beta_2^\vee&=\al_1^\vee+2\al_2^\vee+2\al_3^\vee+3\al_4^\vee+3\al_5^\vee+2\al_6^\vee+\al_7^\vee+\al_8^\vee=\frac 12(1,1,-1,1,1,-1,1,1) \\
\beta_3^\vee&=\al_1^\vee+2\al_2^\vee+2\al_3^\vee+3\al_4^\vee+2\al_5^\vee+2\al_6^\vee+2\al_7^\vee+\al_8^\vee=\frac 12(1,1,1,-1,-1,1,1,1) \\
\beta_4^\vee&=\al_1^\vee+\al_2^\vee+2\al_3^\vee+3\al_4^\vee+3\al_5^\vee+2\al_6^\vee+2\al_7^\vee+\al_8^\vee=\frac 12(-1,-1,-1,1,-1,1,1,1). \\
\end{aligned}
\end{equation}
These are the coroots labelled $89,90,91,92$ in {\it loc. cit.}.

Write $\nu=\sum_{i=1}^8\nu_i\omega_i\in \CC_0^{\frac 12}$, $\nu_i\in \frac 12\bZ_{\ge 0}$ for all $i$. Using (\ref{bound-E8}), we immediately see that necessary conditions are
\[\ 2\nu_2<1,\ 2\nu_3<1,\ 4\nu_4<1,\ 3\nu_5<1,\ 2\nu_6<1,\ 2\nu_7<1,\text{ and } \nu_1+\nu_8<1.
\]
From this, we must have 
\[\nu_2=\nu_3=\nu_4=\nu_5=\nu_6=\nu_7=0,\text{ and }\nu_1,\nu_8\in\{0,\frac 12\}\text{ but not }\nu_1=\nu_8=\frac 12.
\]
This means that there are three remaining points to check:
\begin{itemize}
\item $\nu=0$, this is the origin;
\item $\nu=\frac 12\omega_1$;
\item $\nu=\frac 12\omega_8$;
\end{itemize}
But the last two points are both on the hyperplane $\gamma^\vee=1$, where $\gamma^\vee=\ep_7+\ep_8$ is the highest coroot.

\subsection{$E_7$} We use $\al_i^\vee$, $1\le i\le 7$ from $E_8$, but the coordinates for the fundamental weights are different:
\begin{equation}
\begin{aligned}
\om_1&=(0,0,0,0,0,0,-1,1) &\om_2&=(\frac 12,\frac 12,\frac 12,\frac 12,\frac 12,\frac 12,-1,1), &\om_3&=(-\frac 12,\frac 12,\frac 12,\frac 12,\frac 12,\frac 12,-\frac 32,\frac 32),\\
 \om_4&=(0,0,1,1,1,1,-2,2),& \om_5&=(0,0,0,1,1,1,-\frac 32,\frac 32), &\om_6&=(0,0,0,0,1,1,-1,1),\\
 \om_7&=(0,0,0,0,0,1,-\frac 12,\frac 12).
\end{aligned}
\end{equation}
From \cite[\S7.2.2]{BC-E8}, we see that $r_0=9$ and there are four coroots at level $9$. They are:
\begin{equation}\label{bound-E7}
\begin{aligned}
\beta_1^\vee&=\al_1^\vee+\al_2^\vee+2\al_3^\vee+2\al_4^\vee+2\al_5^\vee+\al_6^\vee=\frac 12(-1,1,-1,1,1,-1,-1,1) \\
\beta_2^\vee&=\al_1^\vee+\al_2^\vee+2\al_3^\vee+2\al_4^\vee+\al_5^\vee+\al_6^\vee+\al_7^\vee=\frac 12(-1,1,1,-1,-1,1,-1,1) \\
\beta_3^\vee&=\al_1^\vee+\al_2^\vee+\al_3^\vee+2\al_4^\vee+2\al_5^\vee+\al_6^\vee+\al_7^\vee=\frac 12(1,-1,-1,1,-1,1,-1,1) \\
\beta_4^\vee&=\al_2^\vee+\al_3^\vee+2\al_4^\vee+2\al_5^\vee+2\al_6^\vee+\al_7^\vee=\ep_5+\ep_6. \\
\end{aligned}
\end{equation}
These are the coroots labelled $46,47,48,49$ in {\it loc. cit.}.

Write $\nu=\sum_{i=1}^7\nu_i\omega_i\in \CC_0^{\frac 12}$, $\nu_i\in \frac 12\bZ_{\ge 0}$ for all $i$. Using (\ref{bound-E7}), we immediately see that necessary conditions are
\[\ 2\nu_3<1,\ 2\nu_4<1,\ 2\nu_5<1,\ 2\nu_6<1, \text{ and } \nu_1+\nu_2+\nu_7<1.
\]
From this, we must have 
\[\nu_3=\nu_4=\nu_5=\nu_6=0,\text{ and }\nu_1,\nu_2,\nu_7\in\{0,\frac 12\}\text{ but only one is nonzero}.
\]

This means that there are four remaining points to check:
\begin{itemize}
\item $\nu=0$, this is the origin;
\item $\nu=\frac 12\omega_7$;
\item $\nu=\frac 12\omega_1$;
\item $\nu=\frac 12\omega_2$;
\end{itemize}
But the last two points are both on the hyperplane $\gamma^\vee=1$, where $\gamma^\vee=-\ep_7+\ep_8$ is the highest coroot. On the other hand, the second point, $\nu=\frac 12\omega_7$ satisfies
\[\langle\gamma^\vee,\frac 12\omega_7\rangle=\frac 12<1,
\]
so it lies in the (open) fundamental alcove, which is always unitary.

\subsection{$E_6$} We use $\al_i^\vee$, $1\le i\le 6$ from $E_8$, but the coordinates for the fundamental weights are different:
\begin{equation}
\begin{aligned}
\om_1&=(0,0,0,0,0,-2/3,-2/3,2/3), &\om_2&=(\frac 12,\frac 12,\frac 12,\frac 12,\frac 12,-\frac 12,-\frac 12,\frac 12),\\
 \om_3&=(-\frac 12,\frac 12,\frac 12,\frac 12,\frac 12,-\frac 56,-\frac 56,\frac 56),
 &\om_4&=(0,0,1,1,1,-1,-1,1),\\ \om_5&=(0,0,0,1,1,-\frac 23,-\frac 23,\frac 23), &\om_6&=(0,0,0,0,1,-\frac 13,-\frac 13,\frac 13).
 \end{aligned}
\end{equation}
Write $\nu=\sum_{i=1}^6\nu_i\omega_i\in \CC_0^{\frac 12}$, $\nu_i\in \frac 12\bZ_{\ge 0}$ for all $i$. Such a parameter is hermitian ($w_0(\nu)=-\nu$) if and only if 
\[\nu_1=\nu_6\text{ and }\nu_3=\nu_5.
\]
From \cite[\S3.5]{Ci-E6}, see also \cite[\S7.2.1]{BC-E8},  we see that $r_0=9$ and there is one coroot at level $9$:
\begin{equation}\label{bound:E6}
\beta^\vee=\al_1^\vee+\al_2^\vee+2\al_3^\vee+2\al_4^\vee+2\al_5^\vee+\al_6^\vee=\frac 12(-1,1,-1,1,1,-1,-1,1).
\end{equation}
From this, it is immediate that necessary conditions are:
\[2\nu_3<1,\ 2\nu_4<1,\ 2\nu_5<1,\text{ and }\nu_1+\nu_2+\nu_6<1.
\]
Hence $\nu_3=\nu_4=\nu_5=0$. Since $\nu_1=\nu_6$, it also follows that $\nu_1=\nu_6=0$. We are left with $\nu_2<1$. The only other point besides the origin is $\nu=\frac 12\omega_2$. But this point lies on the hyperplane $\gamma^\vee=1$, where $\gamma^\vee=\frac 12(1,1,1,1,1,-1,-1,1)$ is the highest coroot.

\subsection{Other fundamental representations} We retain the notation from the previous section. Let $V=V(\lambda^\vee)$ be an irreducible finite-dimensional representation of $G^\vee$ with highest weight $\lambda^\vee\in \ft$.

\begin{lemma}\label{l:half-other}
Let  $\nu\in \CC_0^{\frac 12}\subset \ft^\vee$ be given. Then $\nu$ is half-integral for $V(\lambda^\vee)$ if and only if $\langle\lambda^\vee,\nu\rangle \in \frac 12\bZ$.
\end{lemma}

\begin{proof}
If $\lambda'^\vee$ is a weight of $V(\lambda)$, then we may write $\lambda'^\vee=\lambda^\vee-\sum_{i=1}^\ell a_i\alpha_i^\vee$, where $a_i\in \bZ_{\ge 0}$. Then 
\[\langle \lambda'^\vee,\nu\rangle=\langle\lambda^\vee,\nu\rangle-\sum_{i=1}^\ell a_i\langle\al_i^\vee,\nu\rangle\in \langle\lambda^\vee,\nu\rangle+\frac 12\bZ,
\]
since  $\nu\in \CC_0^{\frac 12}$ and hence by definition $\langle\al_i^\vee,\nu\rangle\in \frac 12\bZ$. The claim follows.
\end{proof}

\begin{prop}\label{p:half-other} Let $G$ and the rest of the notation be as in Theorem \ref{t:half-adjoint}. 
\begin{enumerate} 
\item[(i)] If $G=PGL(n)$, $n$ even, then $\frac 12\omega_{\frac n2}$ is not half-integral for the defining representation $V(\omega_1^\vee)$;
\item[(ii)] If $G=SO(2n+1)$, then $\frac 12 \omega_1$ is not half-integral for the defining representation $V(\omega_n^\vee)$;
\item[(iii)] If $G=PSp(2n)$, then $\frac 12 \omega_n$ is not half-integral for the spin representation $V(\omega_1^\vee)$;
\item[(iv)] If $G=PSO(2n)$, then $\frac 12 \omega_n$ is not half-integral for the spin representation $V(\omega_1^\vee)$, while $\frac 12\omega_1$ and $\frac 12\omega_2$ are not half-integral for the defining representation $V(\omega_n^\vee)$;
\item[(v)] If $G=E_7$, then $\frac 12\omega_7$ is not half-integral for the $56$-dimensional fundamental representation $V(\omega_7^\vee)$.
\end{enumerate}
\end{prop}

\begin{proof} In light of Lemma \ref{l:half-other}, we only need to check the pairing of the listed points with the corresponding highest weight. In the notation from Theorem \ref{t:half-adjoint} and the subsequent subsections, the relevant pairings are as follows:
\begin{enumerate}
\item[(i)] $\langle \omega_1^\vee, \frac 12\omega_{\frac n2}\rangle =\frac 14$;
\item[(ii)] $\langle \omega_n^\vee, \frac 12\omega_1\rangle =\frac 14$;
\item[(iii)] $\langle \omega_1^\vee, \frac 12\omega_n\rangle =\frac 14$;
\item[(iv)] $\langle \omega_1^\vee, \frac 12\omega_n\rangle =\langle \omega_n^\vee, \frac 12\omega_1\rangle =\langle \omega_n^\vee, \frac 12\omega_2\rangle =\frac 14$;
\item[(v)] $\langle \omega_7^\vee, \frac 12\omega_7\rangle=\frac 34$.
\end{enumerate}

\end{proof}

We can now state our conclusions.

\begin{cor}\label{c:half-integral} Let $s=s_c s_\nu$ be a Satake parameter as in (\ref{e:polar}) and let $L(s)$ be the corresponding spherical $G(k)$-representation. Suppose $L(s)$ is generic.
\begin{enumerate}
\item If $L(s)$ is unitary, then $s$ is geometric with respect to the adjoint $G^\vee$-representation and the fundamental $G^\vee$-representations listed in Proposition \ref{p:half-other} if and only if $\nu=0$, i.e., if and only if $L(s)=L(s_c)$ is tempered.
\item If there exists a coroot $\alpha^\vee$ such that $\langle \alpha^\vee,s\rangle=q^m$ for some $m\in \bZ\setminus\{0\}$, then $L(s)$ is not unitary.
\end{enumerate}
\end{cor}

\begin{proof}
For (1), clearly $s$ is geometric with respect to a given $G^\vee$-representation $V$ if and only if $\nu$ is half integral. Without loss of generality, we may assume that $\nu$ is dominant. If $L(s)$ is assumed generic and unitary, we must have $\nu\in \CU_0^g$. By Theorem \ref{t:half-adjoint}, $\nu\in \CU_0^{g,\frac 12}$, a set that was determined explicitly. Using Proposition \ref{p:half-other}, no nonzero point in $ \CU_0^{g,\frac 12}$ can be half-integral for the listed fundamental representations.
The converse is obvious.

Part (2) follows at once from Proposition \ref{p:integral-unitary}.
\end{proof}

\section{Nonsplit unramified groups}\label{s:nonsplit}

In this section, let $G$ be a quasi-simple {\it unramified} group of adjoint type over the nonarchimedean local field $k$. This means \cite[p. 135]{Ca} that
\begin{enumerate}
\item[(a)] $G$ is quasi-split over $k$, and
\item[(b)] there exists an unramified extension $k'$ of $k$ of finite degree $d$ such that $G$ splits over $k'$.
\end{enumerate}
For applications to the automorphic setting, these are the groups that need to be considered. In the previous section, we proved the unitarity results in the case when $G$ is $k$-split. Now, we show how to deduce the analogous results for $G$ unramified from the split case. We follow the approach of \cite[\S6.1-6.2]{BC-ajm}. We note that the assumption that $k$ has characteristic $0$ stated in \cite[\S1.3]{BC-ajm} is unnecessary.

Let $A$ be a maximally split torus of $G(k)$ and $M=C_{G(k)}(A)$, $W(G,A)=N_{G(k)}(A)/M$. Let $K$ be a special maximal compact subgroup attached via the Bruhat-Tits building to $A$. Let $B\supset A$ be a $k$-rational Borel subgroup. Let $X^*(M)$ and $X_*(M)$ denote the lattices of algebraic characters and cocharacters of $M$, respectively. Let $\langle~,~\rangle$ be their
natural pairing. Define the valuation map as before
\begin{equation*}
  \begin{aligned}
\mathsf{val}_M:&~M\to X_*(M),\langle\lambda,\mathsf{val}_M(m)\rangle=\mathsf{val}_k(\lambda(m)),\text{ for all } m\in M,\ \lambda\in X^*(M).    
  \end{aligned}
\end{equation*}
Set $^0\!M=\ker v_M$ and let $\Lambda(M)$ be the image of $\mathsf{val}_M.$ 

The group of unramified characters $\widehat {M}^u$ of $M$ (i.e., characters trivial on $^0\!M$) can be identified with the complex algebraic torus $T'^{\vee}=\text{Spec}~\bC[\Lambda(M)]$. 
For every $s\in T'^{\vee}$, let  $\chi_s\in \widehat M^u$ be the corresponding unramified character and let $X(s)=i_{B(k)}^{G(k)}(\chi_s)$ denote the  unramified principal series. 

As before, denote by $L(s)$ the unique irreducible $K$-spherical subquotient. We have
a  one-to-one correspondence $s\leftrightarrow L(s)$ between $W(G,A)$-orbits in $T'^{\vee}$ and irreducible $K$-spherical representations of $G(k)$.

We need to express this bijection in terms of $^L\!G$.  Let
$\Psi(G^\vee)=(X^*(T^\vee),X_*(T^\vee),\Phi^\vee,\Phi)$ be the root datum of $G^\vee$, and let $W$ be its Weyl group.
The inner class of $G$ defines a  homomorphism $\tau:
Gal(\overline k/k)\to \mathrm{Aut}(\Psi(G))$. Given our assumptions of $G$, we know that the image
$\tau(Gal(\overline k/k))\subset  \mathrm{Aut}(\Psi(G))$ is a cyclic group 
generated by an automorphism of order $d$ ($d\in\{2,3\}$) which, by abuse of notation, we also call
$\tau.$ Fix a choice of 
root vectors in $\mathfrak g^\vee$, $X_{\alpha^\vee},$ for $\alpha\in \Phi.$ The automorphism $\tau$ maps the root space
of $\alpha^\vee$ to the root space of $\tau(\alpha^\vee).$ We can normalize $\tau$ so
that 
$\tau(X_{\alpha^\vee})=X_{\tau(\alpha^\vee)}$, for all $\alpha$.
Hence $\tau$ induces an automorphism of $G^\vee$. 

Two elements $x_1,x_2\in G^\vee$ are called {\it $\tau$-conjugate} if there
exists $g\in G$ such that $x_2=g\cdot x_1\cdot \tau(g^{-1}).$ For a
subset $S\subset G$, denote
\[N_G(S\tau)=\{g\in G:~ g\cdot S\cdot \tau(g^{-1})\subset S\}.\]

 The construction of the L-group is such that we have (see
 \cite[section 6]{Bo}):
\begin{equation}
X^*(T'^{\vee})=X^*(T^\vee)^\tau,\quad W(G,A)=W^\tau.
\end{equation}
In particular, we have an inclusion $X^*(T'^{\vee})\hookrightarrow X^*(T^\vee)$,
which gives a surjection $\eta:T^\vee\to T'^{\vee}.$ 
We recall that by \cite[Lemma 6.4, Proposition 6.7]{Bo}, the map 
\begin{equation}
\eta': T^\vee\rtimes\langle\tau\rangle\to T'^{\vee},\quad \eta'((s,\tau))=\eta(s),
\end{equation}
induces a bijection of $(N(T^\vee\tau),\tau)$-conjugacy classes of elements
in $T^\vee$ (hence of semisimple $\tau$-conjugacy classes in $G^\vee$) and $W^\tau$-conjugacy classes of elements in $T'^{\vee}.$

Let
\begin{equation}G^\vee(\tau)=\{g\in G^\vee: \tau(g)=g\}.\end{equation} 

This is a connected reductive group (recall that $G^\vee$ is assumed to be simply connected).  Let  $G(\tau)$ denote the split $k$-adic group whose Langlands dual  is $G^\vee(\tau)$. The explicit cases are listed in Table \ref{t:1}.

\begin{table}[h]\caption{Nonsplit quasisplit unramified $G$ and corresponding split $G(\tau)$.\label{t:1}}
\begin{tabular}{|c|c|c|c|}
\hline
Type of $G$ &Label in \cite{Ti} &Order of $\tau$ &Type of $G(\tau)$\\
\hline
$A_{2n-1}$ &$^2\!A_{2n-1}'$ &$2$ &$B_n$\\
\hline
$A_{2n}$ &$^2\!A_{2n}'$ &$2$ &$C_n$\\
\hline
$D_n$ &$^2\!D_n$ &$2$ &$C_{n-1}$\\
\hline
$D_4$ &$^3\!D_4$ &$3$ &$G_2$\\
\hline
$E_6$ &$^2\!E_6$ &$2$ &$F_4$\\
\hline
\end{tabular}
\end{table}

As in section \ref{classificationsection}, for every $\nu\in \mathfrak t_{\mathbb R}'^\vee=(\mathfrak t_{\mathbb R}^\vee)^\tau$, let $X^G(\nu)=X({\exp(\nu\log(q)})$ be the unramified principal series for $G(k)$. Let also $X^{G(\tau)}(\nu)$ be the unramified principal series for $G(\tau)(k)$. Let $L^G(\nu)$ and $L^{G(\tau)}(\nu)$ be the corresponding spherical subquotients. The following statement is a particular case of \cite[Corollary 6.2.6]{BC-ajm}

\begin{thm}\label{t:quasi-corr} For every $\nu\in \mathfrak t_{\mathbb R}'^\vee$, the irreducible representation $L^G(\nu)$ is (generic) unitary if and only if $L^{G(\tau)}(\nu)$ is (generic) unitary.
\end{thm}

From this we deduce immediately the following analogue of Corollary \ref{c:half-integral}.

\begin{cor}\label{c:unramified}
Suppose $G$ is a quasi-simple unramified nonsplit adjoint $k$-group and that $G$ is not the quasisplit form of $PSU(2n+1)$. Let $s=s_c s_\nu\in T'^\vee$ be a Satake parameter and $L^G(s)$ the corresponding spherical $G(k)$-representation. Suppose $L^G(s)$ is generic and unitary. Then $s$ is geometric with respect to all finite-dimensional representations of $G^\vee\rtimes\langle \tau\rangle$ if and only if $\nu=0$, i.e., if and only if $L^G(s)=L^G(s_c)$ is tempered.
\end{cor}

\begin{proof}
The generic spherical $G(k)$-representation $L^G(\nu)$ is unitary, hence by Theorem \ref{t:quasi-corr}, $L^{G(\tau)}(\nu)$ is generic unitary. We wish to apply Theorem \ref{t:half-adjoint} and Proposition \ref{p:half-other} to $L^{G(\tau)}(\nu)$. For this, we need to exhibit finite-dimensional representations of $G^\vee\rtimes \langle\tau\rangle$ such that their restrictions to $G^\vee(\tau)$ contain the adjoint representation of $G^\vee(\tau)$ and the representations listed in Proposition \ref{p:half-other}. 

Firstly, the adjoint  $G^\vee$-representation is $\tau$-stable and contains in its restriction the adjoint representation of $G^\vee(\tau)$.  Next, it remains to discuss the classical groups listed in Proposition \ref{p:half-other}; in all cases $\tau^2=1$. 

\begin{enumerate}
\item $G^\vee=SL(2n,\bC)$, $G^\vee(\tau)=Sp(2n,\bC)$, $G(\tau)=SO(2n+1,k)$. The vector $G^\vee$-representation $V=\bC^{2n}$  restricts to the defining representation $V(\omega_n^\vee)$ of $G^\vee(\tau)$.
\item $G^\vee=Spin(2n,\bC)$, $G^\vee(\tau)=Spin(2n-1,\bC)$, $G(\tau)=PSp(2n-2,k)$. The two spin $G^\vee$-representations restrict both to the spin representation of $G^\vee(\tau)$.  
\end{enumerate} 
In this way, Proposition \ref{p:half-other} can be applied to $G^\vee(\tau)$ and it follows that $\nu=0$.

\end{proof}

\begin{remark}\label{PSU}
The methods in this section are insufficient to establish the claim in Corollary \ref{c:unramified} in the case when $G$ is the quasisplit form of $PSU(2n+1)$. In that case, $G^\vee=SL(2n+1,\bC)$, $G^\vee(\tau)=SO(2n+1,\bC)$, and $G(\tau)=Sp(2n,k)$. The adjoint representation of $G^\vee$ contains in its restriction the adjoint representation of $G^\vee(\tau)$, which means that the argument in the proof of the corollary implies that the only possible geometric Satake parameters for $G(k)$ have $\nu=0$ or $\nu=\frac 12\omega_n$, where $\omega_n$ is the $G(\tau)$-weight (the notation as in subsection \ref{sub:psp}). But in the case of $PSp(2n,k)$, the point $\nu=\frac 12\omega_n$ is ruled out by the spin representation of the dual group $Spin(2n+1,\bC)$, see Proposition \ref{p:half-other}(iii), and the spin representation does not descend to $SO(2n+1,\bC)$ (in particular, it can't appear in the restriction of a finite-dimensional $SL(2n+1,\bC)$-representation). 
\end{remark}

\section{Complementary series in the spherical unitary spectrum of a $p$-adic group}

Let $\mathsf k$ denote a nonarchimedean local field with finite residue field $\mathbb F_q$. There is no assumption on the characteristic of $\mathsf k$ in this section, but our main application will be the specialization $\mathsf k=K_w$, in terms of the notation from the previous sections. Let $G$ be a quasi-simple $\mathsf k$-split group. Let $T$ be a maximally split torus of $G(\mathsf k)$ and $W=N_{G(\mathsf k)}(T)/T$ the Weyl group. Let $K$ be a maximal hyperspecial compact open subgroup of $G(\mathsf k)$ attached via the Bruhat-Tits building to $T$. As in the beginning of \S\ref{classificationsection}, for every $s\in T^{\vee}$, let  $\chi_s$ be the corresponding unramified character of $T$ and let $X(s)=i_{B(\mathsf k)}^{G(\mathsf k)}(\chi_s)$ denote the  unramified principal series. 
Recall that $L(s)$ the unique irreducible $K$-spherical subquotient. The Satake correspondence states that there is a 
one-to-one correspondence 
$s\longleftrightarrow L(s)$ between $W$-orbits in $T^\vee$ and irreducible $K$-spherical representations.

\subsection{Nilpotent orbits} Let $I\subset K$ be an Iwahori subgroup of $G(\mathsf k)$. As recalled in the proof of Corollary \ref{locgeneric}, the smooth irreducible $I$-spherical representations of $G(\mathsf k)$ are in one-to-one correspondence with $G^\vee$-conjugacy classes of Kazhdan-Lusztig triples $(s,e^\vee,\rho)$, where
\begin{itemize}
\item $s\in G^\vee$ is semisimple;
\item $e^\vee\in \mathfrak g^\vee$ is nilpotent and $\mathrm{Ad}(s) e^\vee=q e^\vee$;
\item $\rho$ is an irreducible representation of the group of components of the centralizer of $s, e^\vee$ in $G^\vee$ of Springer type, see \cite[Theorem 7.12]{KL} for the precise statement. 
\end{itemize}
In the case when $G$ has connected center, this classification is proved in \cite{KL}. The isogeny assumption was removed in \cite[Theorem 1]{Re-isogeny}.

Let $\pi(s,e^\vee,\rho)$ denote the irreducible representation parametrized by $(s,e^\vee,\rho)$. In this correspondence, the $K$-spherical representations correspond precisely to $e^\vee=0$ and $\rho$ trivial, i.e., $L(s)=\pi(s,0,\mathsf{triv})$. The generic representations are $\pi(s,e^\vee,\rho)$ such that $\rho$ is trivial and $e^\vee$ lies in the unique open $G^\vee(s)$-orbit in $\mathfrak g^\vee_q$, the $q$-eigenspace of $\mathrm{Ad}(s)$. The Aubert-Zelevinsky involution preserves the semisimple parameter $s$ and maps $K$-spherical representations to generic representations, that is
\[AZ(L(s))=\pi(s,e^\vee,\mathsf{triv}),
\]
where $e^\vee$ lies in the open  $G^\vee(s)$-orbit in $\mathfrak g^\vee_q$. For the organization of the unitary dual, it is useful to introduce the following notation. For every nilpotent $G^\vee$-orbit $\CO^\vee$, set
\begin{equation}
\mathsf{Sph}(G,\CO^\vee)=G^\vee\text{-conjugacy classes }\{s\in G^\vee\mid AZ(L(s))=\pi(s,e^\vee,\mathsf{triv}),\text{ where }e^\vee\in \CO^\vee\}.
\end{equation}
Notice that the sets $\mathsf{Sph}(G,\CO^\vee)$ as $\CO^\vee$ varies over the finite set of nilpotent orbits of $\mathfrak g^\vee$ give a partition of the $K$-spherical dual of $G(\mathsf k)$. Fix a  Lie triple $(e^\vee,h^\vee,f^\vee)$ for $\CO^
\vee$. Then $\mathsf{Sph}(G,\CO^\vee)$ contains a distinguished parameter, which we call the \emph{central point}, this is
\[s=q^{h^\vee/2}\in \mathsf{Sph}(G,\CO^\vee).
\]
The corresponding spherical representation is the $AZ$-dual of a tempered representation and it is unitary.

Let $T^\vee=T^\vee_{\bR} T^\vee_c$ denote the polar decomposition of $T^\vee$. 
If $s\in T^\vee_{\bR}$, let $\Re(s)$ denote the unique element of $\mathfrak t^\vee_{\bR}=\Lie( T^\vee_{\bR})$ such that $s= q^{\Re(s)} \cdot s_c$, for some $s_c\in T^\vee_c$. Denote
\[\mathsf{Sph}(G,\CO^\vee)_0=G^\vee\text{-conjugacy classes }\{\Re(s)\mid s\in \mathsf{Sph}(G,\CO^\vee)\}.
\]

\subsection{Complementary series}\label{s:CS} Let $\mathcal{SU}(G)$ denote the spherical unitary dual of $G(\mathsf k)$, the set of Satake parameters $s$ such that the $K$-spherical representation $L(s)$ is unitary. Set
\[\mathsf{CS}(G,\CO^\vee)=\mathcal{SU}(G)\cap \mathsf{Sph}(G,\CO^\vee).
\]
By definition $\mathcal{SU}(G)=\sqcup_{\CO^\vee} \mathsf{CS}(G,\CO^\vee)$. We call $\mathsf{CS}(G,\CO^\vee)$ the \emph{spherical complementary series} of $G(\mathsf k)$ attached to $\CO^\vee$. Denote also
\[\mathsf{CS}(G,\CO^\vee)_0=G^\vee\text{-conjugacy classes }\{\Re(s)\mid s\in \mathsf{CS}(G,\CO^\vee)\},
\]
the real complementary series attached to $\CO^\vee$.
Denote $Z(\CO^\vee)=C_{G^\vee}(e^\vee,h^\vee,f^\vee)$, this is a possibly disconnected reductive group with Lie algebra $\mathfrak z(\CO^\vee)$. Let $G(\CO^\vee)$ denote the split $\mathsf k$-group which is dual to the identity component $Z(\CO^\vee)^0$. 

Following \cite[\S1.3,\S2.5]{BC-E8}, let $\mathfrak a^\vee$ be a Cartan subalgebra of $\mathfrak z(\CO^\vee)$ such that $\mathfrak a^\vee\subset \mathfrak t^\vee$. If $\Re(s)\in \mathsf{Sph}(G,\CO^\vee)_0$, we may arrange that $\Re(s)\in \mathfrak t^\vee$ (up to conjugacy by the Weyl group $W$) and we write
\begin{equation}
\Re(s)=\frac 12 h^\vee+\nu\text{ for some }\nu\in \mathfrak a_\bR^\vee.
\end{equation}
We regard $\nu$, as we may, as a parameter in $\mathsf{Sph}(G(\CO^\vee),0)_0$. The classification of the spherical unitary dual from \cite{Ba,BC-E8,Ci-F4,Ci-E6} implies in particular, that with one exception, $\mathsf{CS}(G,\CO^\vee)_0\subseteq \mathsf{CS}(G(\CO^\vee),0)$. More precisely:

\begin{thm}\label{t:BC}
If $\CO^\vee$ is not the nilpotent orbit $4A_1$ in $E_8$, and $s\in \mathsf{CS}(G,\mathcal O^\vee)$ has $\Re(s)=\frac 12 h^\vee+\nu$, with $\nu\in \mathsf{Sph}(G(\CO^\vee),0)_0$, then $\nu\in \mathsf{CS}(G(\CO^\vee),0)_0$.
\end{thm} 
In the case $4A_1$ in $E_8$, the complementary series has an extra unitary region compared to that for the centralizer, this will be treated separately in subsection \ref{s:4A1}.

\begin{defn}\label{d:extraneous}   When $s = \varphi(\Frob)$ is geometric and $L(s)$ is unitary, but $\Re(s)$ is not the central point, we say $\Re(s)$ is an {\it extraneous point}. \end{defn}

With the description in Theorem \ref{t:BC}, we can prove the main result of the section.

\begin{thm}\label{t:central}
Assume $G$ has connecter center. Suppose $s=\varphi(\Frob)\in \mathsf{CS}(G,\CO^\vee)$ is geometric in the sense of Definition \ref{algebraic}. The extraneous points in the sense of Definition \ref{d:extraneous} occur in the following cases:
\begin{enumerate}
\item $G^\vee=Sp(2n)$, $\CO^\vee=(\underbrace{1,\dots,1}_{r_1},\underbrace{2,\dots,2}_{r_2},\dots)$, $\sum_i i r_i=2n$, $r_i$ even if $i$ odd, $Z(\CO^\vee)=\prod_{i\text{ odd}} Sp(r_i)\times \prod_{j\text{ even}} O(r_j)$,
\[\Re(s)=\frac 12 h^\vee+ \frac 12\sum_{j\text{ even}, r_j\ge 3} \epsilon_j \omega_{\lfloor r_j/2\rfloor}^j,\quad \epsilon_j\in\{0,1\},
\]
where not all $\epsilon_j$ are zero, and $\omega_{\lfloor r_j/2\rfloor}^j$ is the fundamental coweight of $so(r_j)$ as in Proposition \ref{p:half-other}.

\item $G^\vee=Spin(n)$, $\CO^\vee=(\underbrace{1,\dots,1}_{r_1},\underbrace{2,\dots,2}_{r_2},\dots)$, $\sum_i i r_i=n$, $r_i$ even if $i$ even, $\mathfrak z(\CO^\vee)=\oplus_{i\text{ odd}} \mathfrak{so}(r_i)\bigoplus \oplus_{j\text{ even}} \mathfrak{sp}(r_j)$,
\[\Re(s)=\frac 12 h^\vee+ \frac 12\sum_{i\text{ odd}, r_i\ge 3} \epsilon_i \omega_{\lfloor r_i/2\rfloor}^i,\quad \epsilon_i\in\{0,1\},
\]
where $\sum_{i\text{ odd}, r_i\ge 3} \epsilon_i \equiv 0$ (mod $2$), not all $\epsilon_i$ are zero, and $\omega_{\lfloor r_i/2\rfloor}^i$ is as in (1).

\smallskip

\item $G^\vee=F_4$

\begin{tabular}{|c|c|c|}
\hline
$\CO^\vee$ &$\mathfrak z(\CO^\vee)$ &$\Re(s)$\\
\hline
$B_3$ &$A_1$ & $\omega_1+\omega_2+\frac 12\omega_3+\frac 12\omega_4$\\
\hline
$A_1+\widetilde A_1$ &$A_1+A_1$ &$\frac 12\omega_1+\frac 12\omega_4$\\
\hline
\end{tabular}

\smallskip

\item $G^\vee=E_7$, $\CO^\vee=D_5(a_1)+A_1$, $\mathfrak z(\CO^\vee)=A_1$,
\[\Re(s)=\frac 12 \omega_1+\frac 12\omega_2+\frac 12\omega_3+\frac 12\omega_6+\frac 12\omega_7.
\]

\smallskip

\item $G^\vee=E_8$

\begin{tabular}{|c|c|c|}
\hline
$\CO^\vee$ &$\mathfrak z(\CO^\vee)$ &$\Re(s)$\\
\hline
 $D_6$ &$B_2$ &$\omega_1+\frac 12\omega_4+\frac 12\omega_6+\omega_8$\\
 \hline
$A_6$ &$2A_1$ &$\frac 12 \omega_1+\frac 12\omega_4+\frac 12\omega_7+\frac 12\omega_8$\\
\hline
$A_4+A_2$ &$2A_1$ &$\frac 12\omega_1+\frac 12\omega_5+\frac 12\omega_8$\\
\hline
$A_2+2A_1$ &$B_3+A_1$ &$\frac 12\omega_2+\frac 12\omega_8$\\
\hline
\end{tabular}

\end{enumerate}
\end{thm}

\begin{remark}
The smallest example in the previous theorem is $G^\vee=Sp(6)$, $\CO^\vee=(2,2,2)$, where $\frac 12h^\vee=(\frac 12,\frac 12,\frac 12)$ in the usual coordinates, and the extraneous point is $s=(0,\frac 12,1)$.
\end{remark}

\begin{proof}
If $\CO^\vee=4A_1$ in $E_8$, the proof is given in subsection \ref{s:4A1}. Otherwise, write $\Re(s)=\frac 12 h^\vee+\nu$, where $\nu\in \mathsf{CS}(G(\CO^\vee),0)_0$. Notice that, by $\mathfrak{sl}(2)$-theory, $\frac 12 h^\vee$ is half-integral with respect to any (finite-dimensional) representation $(V,\rho)$ of $\mathfrak g^\vee$. Since  $s$ is geometric for the group $G$, it follows that $\nu$ is half-integral for every $(V,\rho)$ of $\mathfrak g^\vee$. In particular, $\nu$ is half-integral for the adjoint representation of $\mathfrak g$, and since $\mathfrak z(\CO^\vee)\subset \mathfrak g^\vee$, $\nu \in \mathsf{CS}(G(\CO^\vee),0)_0$ is half-integral for the adjoint representation of $\mathfrak z(\CO^\vee)$ as well. Applying Theorem \ref{t:half-adjoint} to $\mathfrak z(\CO^\vee)$, it follows that if $\nu\neq 0$, then $\nu$ must equal one of the exceptions listed in Theorem \ref{t:half-adjoint} on each of the simple factors of $\mathfrak z(\CO^\vee)$. 

To rule out the remaining nonzero points $\nu$, we need to proceed as in the proof of Proposition \ref{p:half-other}. The difficulty comes from the fact that a priori there is no reason for all of the fundamental representations of $\mathfrak z(\CO^\vee)$ to occur in the restriction of some finite-dimensional representation of $\mathfrak g^\vee$. In fact this is not true, for example $\mathfrak z(\CO^\vee)$ could be $\mathfrak{so}(m)$ inside $\mathfrak g^\vee=\mathfrak{sp}(2m)$, and the spin module(s) of $\mathfrak{so}(m)$ do not appear in the restriction of any $\mathfrak{sp}(2m)$ finite-dimensional module. Therefore, we do a case-by-case analysis.

\smallskip

(1) If $G^\vee=GL(n)$ and $\CO^\vee=(\underbrace{1,\dots,1}_{r_1},\underbrace{2,\dots,2}_{r_2},\dots)$, $\sum_i i r_i=n$, the centraliser $Z(\CO^\vee)=\prod_i GL(r_i)^\Delta_i$, where $GL(r_i)^\Delta_i$ means a copy of $GL(r_i)$ embedded diagonally in the direct product $GL(r_i)^i$. The standard representation of $GL(r_i)$ occurs in the restriction of the standard representation of $GL(n)$, hence Proposition \ref{p:half-other}(i) can be apply to conclude that $\nu=0$.

\smallskip

(2) If $G^\vee=Sp(2n)$ and $\CO^\vee$ as in the statement of the theorem, write $\nu=\sum_{i\text{ odd}} \nu_i+\sum_{j\text{ even}}\nu_j$ corresponding to the simple factors of the centraliser. We can use the standard representation of $G^\vee$ to conclude that $\nu_i=0$ for all $i$ odd. The remaining nonzero points are $\nu_j=\frac 12 \omega^j_{\lfloor r_j/2\rfloor}$, when $j$ is odd and $r_j\ge 3$, and all of these points are half-integral for $G^\vee$.

\smallskip

(3) If $G^\vee=Spin(n)$, we proceed as in (2). Using the standard representation of $G^\vee$, we find that $\nu_j=0$ for all $j$ even. Using the spin module(s) of $G^\vee$, we rule out the case $\sum_{i\text{ odd}, r_i\ge 3} \epsilon_i\equiv 1$ (mod $2$). The remaining points are all half integral for $G^\vee$.

\smallskip

(4) For the exceptional groups, we list all the orbits and their corresponding centralisers (see the tables in \cite{Ca}). For each simple factor in the centraliser, we list the possible nonzero half-integral points as given by Proposition \ref{p:half-other}. This way we obtain a list of points for each nilpotent orbit which are not ruled out by the adjoint representation of the centraliser. Suppose $\nu$ is such a point. We then compute the values $\alpha^\vee(\nu)$, where $\alpha^\vee$ are the simple roots of $\mathfrak g^\vee$. If any of these values is not in $\frac 12\mathbb Z$, the point $\nu$ is ruled out by the adjoint representation of $G^\vee$. Otherwise, we compute $\omega^\vee(\nu)$ for every fundamental weight $\omega^\vee$ of $G^\vee$. If any of these values are not in $\frac 12\mathbb Z$, again $\nu$ is ruled out. The remaining points are the ones listed in the theorem: two points for $F_4$, one for $E_7$, and four for $E_8$.

\end{proof}

\begin{remark}  Theorem \ref{t:central} gives a complete list of orbits $\CO^\vee$ that admit extraneous points.  For most $\CO^\vee$ there are no such points.  For example, $E_8$ has $70$ nilpotent orbits, but only the four listed in the Theorem admit extraneous points; there are no extraneous points for  $G_2$ or $E_6$.  Our main result below is not optimal in that it does not rule out extraneous points altogether, as would be expected from the Arthur Conjectures; however, it does say that if one unramified component of a cuspidal automorphic representation has extraneous Satake parameter, then it is the same for every unramified component.
\end{remark}

\subsection{The orbit $\CO^\vee=4A_1$ in $E_8$}\label{s:4A1} In this case the $Z(\CO^\vee)$ is of type $C_4$. In coordinates, see \cite[Table 13 and \S6.4.1]{BC-E8}, $\frac 12h^\vee=(0,1,-\frac 12,\frac 12,-\frac 12,\frac 12,0,0)$ and a real Satake parameter can be written as $s=\frac 12h^\vee+(0,0,\nu_1,\nu_1,\nu_2,\nu_2,-\nu_3+\nu_4,\nu_3+\nu_4)$, so that
\[\mathsf{Sph}(E_8,4A_1)=\{0\le \nu_1\le \nu_2\le \nu_3\le\nu_4\}.
\]
The complementary series is
\begin{equation}\mathsf{CS}(E_8,4A_1)=\{\nu_4<\frac 12\}\cup\{\nu_1+\nu_4<1,\nu_2+\nu_3<1,\nu_2+\nu_4>1,-\nu_1+\nu_3+\nu_4<\frac 32<\nu_1+\nu_3+\nu_4\},
\end{equation}
while $\mathsf{CS}(C_4,0)=\{\nu_4<\frac 12\}$. We need to analyze the extra unitary region in $\mathsf{CS}(E_8,4A_1)$.

The half-integral condition for the adjoint representation implies that $\nu_1,\nu_2,\nu_3,\nu_4\in \frac 14 \bZ_{\ge 0}$ are all either in $\frac 12 \bZ_{\ge 0}$ or are all in $\frac 12\bZ+\frac 14$. In the first case, $\nu_i\in \bZ+\frac 12$, and since $\nu_1+\nu_4<1$, we must have $\nu_4<1$ and hence $\nu_4=\frac 12$. But $\nu_2+\nu_4>1$ which implies $\nu_4>\frac 12$, contradiction.

In the second case, $\nu_i\in \bZ+\frac 14$ and using $\nu_2+\nu_4>1$, we see again $\nu_4>\frac 12$. On the other hand, $\nu_1+\nu_4<1$ gives $\nu_4<1$. So the only possibility is $\nu_4=\frac 34$. Since $\nu_2+\nu_4>1$, we have $\nu_2>\frac 14$, which means $\nu_2=\frac 34$, but then also $\nu_3=\frac 34$ and this contradicts $\nu_2+\nu_3<1$. 

In conclusion, the extra unitary region does not have any half-integral parameters. For the region $\{\nu_4<\frac 12\}$ that comes from $\mathsf{CS}(C_4,0)$, the same argument as in the proof of Theorem \ref{t:central} for $\CO^\vee\neq 4A_1$ applies.

\section{Parameters of discrete automorphic representations}

\subsection{Optimistic formulations of the Galois parametrization of discrete automorphic representations}
Notation is as in previous sections.  
We let 
$$\CA_{disc}(G,U,E) = \CA_{disc}\cap \CA(G,U,E) \subset \CA(G,U,E)$$
where $\CA_{disc}$ is the space of square-integrable automorphic forms.  The following Conjecture is the analogue for the discrete spectrum of Proposition \ref{cuspfinite}.

 \begin{conj}\label{discfinite}  Suppose $E$ is a noetherian $\ZZ[\frac{1}{p}]$ algebra with a fixed embedding in the algebraically
 closed field $C$.  
 Then $\CA_{disc}(G,U,E)$ is a finite $E$-module.

 Moreover, $\CA_{disc}(G,U,E)\otimes C \isoarrow \CA_{disc}(G,U,C).$
  \end{conj}
  
When $E = \CC$ the finiteness can probably be derived from the constructions in \cite{MW}, but to our knowledge this has never been written down explicitly.

  
Assume Conjecture \ref{discfinite}; then for any field $E$ of characteristic zero, the action of $\CH^S(G,U)$ on $\CA_{disc}(G,U,E)$ is semisimple.  The proof of Lemma \ref{algHecke} carries over and we have

\begin{lemma}\label{algHeckedisc}
\begin{enumerate}
\item The eigenvalues of $\CH^S(G,U)$ on $\CA_{disc}(G,U,C)$ are algebraic numbers for any algebraically closed $C$.
\item Any  irreducible discrete automorphic representation of $G(\ad_K)$ has a model over a number field.
\end{enumerate}
In particular, 
$$\CA_{disc}(G,\bar{\QQ}) \isoarrow \bigoplus_{\pi}  m(\pi) \pi;   \CA_{disc}(G,U,\bar{\QQ}) \isoarrow \bigoplus_{\pi} m_U(\pi) \pi^U$$
where $\pi$ runs over the irreducible $\bar{\QQ}$-representations of $G(\ad_K)$ and $m(\pi) \geq 0$ is an integer multiplicity, and we make the convention
that $m_U(\pi) = 0$ if $\pi^U = 0$.  
\end{lemma}

At this point we would like to conjecture that, to any irreducible $\pi \subset \CA_{disc}(G,\bar{\QQ})$ -- in other words, with $m(\pi) \geq 0$ -- and for any prime $\ell \neq p$, we can assign a global $\ell$-adic semi-simple Langlands parameter
$$\sigma_\pi = \sigma_{\ell,\pi}:  Gal(K^{sep}/K) \ra \LGr(\Qlb),$$
well defined up to conjugacy by the dual group.
However, this is not even correct for cuspidal $\pi$ when $m(\pi) > 1$.  Lafforgue's construction in that situation does indeed provide a collection of {\it local} Galois parameters 
$$\{\sigma_{\pi,w}: Gal(K_w^{sep}/K_w) \ra \LGr(\Qlb)\}$$
at all places $w$ of $K$.  Our construction, however, requires that the $\sigma_{\pi,w}$ be obtained (up to conjugacy) as the restriction of a global $\sigma_\pi$ to a decomposition group at $w$.

We have seen in \S \ref{sec_param} that, in order to obtain the global parameter, Lafforgue needs to decompose the Hecke eigenspace corresponding to $\pi$ under the (semisimple) action of the algebra $\CB(G,U,\Qlb)$ of excursion operators.  
For the moment we assume there is such an action, and that it gives rise to a decomposition
\begin{equation}\label{paramdisc}
\CA_{disc}(G,U,\Qlb) \isoarrow \bigoplus_{\sigma} \CA_{disc}(G,U,\Qlb)_{\lambda(\sigma)},
\end{equation}
where $\lambda(\sigma)$ runs over characters of $\CB(G,U,\Qlb)$ that occur non-trivially in $\CA_{disc}(G,U,\Qlb)$, and $\sigma$ designates a semisimple Langlands parameter as before.   As in \S \ref{sec_param}, we let $\lambda_\CH(\sigma)$ denote the restriction of $\lambda(\sigma)$ to the $\Qlb$-Hecke algebra, and we attach 
a Galois parameter $\sigma_\pi$ to a given $\pi$ with the following Definition:
\begin{defn}\label{belongs}  Let $\pi \subset \CA_{disc}(G,\bar{\QQ})$ with $m_U(\pi) \neq 0$.  We say $\pi$ {\bf belongs to} the parameter $\sigma$ if $\pi^U$ is contained in $\CA_{disc}(G,U,\Qlb)_{\lambda(\sigma)}$; then 
$\lambda_\CH(\sigma)$ coincides with the action of $\CH^S(G,U,\bar{\QQ})$ on $\pi^U$.  
\end{defn}

Let $\Phi_{disc}(U)$ denote the set of $\sigma$ such that $\CA_{disc}(G,U,\Qlb)_{\lambda(\sigma)} \neq 0$.    Let $\pi$ be a discrete automorphic representation that belongs to $\sigma$.  For $w \notin S(U)$ we let  $\lambda_w(\sigma)$ denote the restriction of $\lambda(\sigma)$ to $\CH_w$.  
We let $\LGr^{ss}(\Qlb)$ denote the set of semisimple elements of $\LGr(\Qlb)$, $[\LGr^{ss}(\Qlb)]$ the set of semisimple conjugacy classes in $\LGr(\Qlb)$.
For $\sigma \in \Phi(U)$ and $w \notin S(U)$, let 
$$\alpha_w(\sigma) = [\sigma(\Frob_w)] \in [\LGr^{ss}(\Qlb)],$$
where the brackets around $\sigma(\Frob_w)$ denote the conjugacy class of the image under $\sigma$ of a choice of Frobenius element at $w$.  The following Conjecture is known for cuspidal representations.

\begin{conj}[Lafforgue]\label{SeqT2}  
\begin{enumerate}  
\item There is a decomposition as in \eqref{paramdisc}.

\item Fix $\sigma \in \Phi(U)$.  
Let
$s_{w,\ell}(\sigma) \in [\LGr^{ss}(\Qlb)]$ be the Satake parameter corresponding to $\lambda_w(\sigma)$.  Then 
$$s_{w,\ell}(\sigma) = \iota_\ell(s_w(\sigma))$$
for an $s_w(\sigma) \in [\LGr^{ss}(\bar{\QQ}]$
and
$$\alpha_w(\sigma) = s_w(\sigma).$$
\end{enumerate}
\end{conj}

\subsection{Review of results of Lafforgue and Xue}

Let $\CA_c(G,U,\Qlb)$ denote the compactly supported $\Qlb$-valued functions on $X_U$.  This is a module over $\CH^S(G,U)$.  
We have already noted that the Conjectures of the previous section were proved by Lafforgue for cuspidal automorphic representations.  This was extended by C. Xue to quotients of $\CA_c(G,U,\Qlb)$ by ideals of finite codimension in $\CH^S(G,U)$ \cite{X}.  Some non-cuspidal discrete automorphic representations $\pi$ arise in such quotients, and then the Conjectures of the previous section are valid for such $\pi$.  The parametrization hypothesis in Theorem \ref{mainthmArthur} applies to cuspidal representations and for the non-cuspidal representations to which Xue's construction applies.  However, the arguments are valid for any $\pi$ for which the Conjectures above are valid.

\section{Proof of Theorem \ref{mainthmArthur}}

The proof is essentially identical to the proof in the generic case.  Recall we are restricting our attention to the case of split groups.  We pick up the proof of Theorem \ref{mainthm} after \eqref{qnumbers}, which we reproduce here:  for any
 weight $\alpha$ of $\tau$, we have
 \begin{equation}\label{qnumbers2} \langle \alpha,\nu_w(\sigma) \rangle \in \frac{1}{2}\ZZ.
 \end{equation}
When $\pi_v$ is not generic, however, but is only assumed to be attached to the nilpotent orbit $\CO^\vee$ containing $N$, we cannot apply the classification of the generic unitary spherical spectrum in \S \ref{classificationsection}.  Thus it is not generally the case that $\nu_v(\sigma) = 0$, and thus the irreducible components of the semi-simple local system $\CL(\pi,\tau)$, with $\tau$ as in the proof of Theorem \ref{mainthm}, are generally not $\iota$-pure of weight $0$.   Instead, as in Proposition \ref{globalweights}, we let $n(i,\pi,\tau)$ denote the rank of the weight $i$ component $\CL(\pi,\tau)$.   We thus need to look more closely at the possible weights that can arise.

\subsection{Weight patterns}  Fix a nilpotent orbit $\CO^\vee$ and let $s \in \mathsf{CS}(G,\CO^\vee)$.    Suppose $s$ is geometric in the sense of Definition \ref{algebraic}.  We define the {\it weight pattern} of $s$ to be the set of functions indexed by irreducible representations $\tau$ of $G^\vee$
$$n_{\tau,s}:  \ZZ \ra \ZZ_{\geq 0};  ~~~ n_{\tau,s}(i) = \dim gr_i L(\varphi,\tau)$$
for any unramified Weil-Deligne parameter $\varphi$ for the local field $F$ such that $\varphi(\Frob) = s$. 

\begin{prop}\label{p:weight-pattern}  The weight pattern of a geometric $s$ determines $\Re(s)$ (up to conjugation).
\end{prop}
\begin{proof}
Since the weight pattern is determined by $\Re(s)$, we may assume without loss of generality that $s=\Re(s)\in T^\vee_{\mathbb R}$. Suppose $s'\in T^\vee_{\mathbb R}$ is another semisimple element with the same weight pattern as $s$. If $V$ is a finite-dimensional representation of $G^\vee$ with character $\chi_V$, then 
\[\chi_V(s)=\sum_{i\in\bZ} n_{\tau,s}(i) q^{i/2},
\]
hence $\chi_V(s)=\chi_V(s')$, for all finite-dimensional $G^\vee$-representations $V$. Via Weyl's character formula, the ring of Weyl group invariant functions $\CO(T^\vee)^W$ is generated by the set of characters $\chi_V$. Therefore, for every $f\in \CO(T^\vee)^W$, $f(s)=f(s')$. Since $\CO(T^\vee)^W$ separates the $W$-conjugacy classes in $T^\vee$, it follows that $s,s'$ are $W$-conjugate.
\end{proof}

Now we conclude the proof of Theorem \ref{mainthmArthur} as in \S \ref{mtgc}.    First note that, just as in the proof of Theorem \ref{mainthm},  we can reduce to the case where $G$ is adjoint by restricting Galois parameters to the Galois group of a finite extension $K'/K$.  Thus the list of extraneous points to which the final sentence of Theorem \ref{mainthmArthur} refers is the list given in Theorem \ref{t:central}.  

The hypothesis that $\pi_u$ is geometric shows that we are in the situation of Proposition \ref{globalweights}.   Thus if $w$ is any unramified place for $\pi$, $\varphi_w$ is the semisimple parameter of $\pi_w$, and $\tau$ is an irreducible representation of ${}^LG$,  the dimension $n(i,\pi,\tau,w)$ of the weight $i$ space for $\tau\circ\varphi_w(\Frob_w)$, with $\Frob_w$ a Frobenius element at $w$, is necessarily $n(i,\pi,\tau)$.  Since this is independent of the unramified place $w$, we must have $n(i,\pi,\tau,w) = n(i,\pi,\tau,v)$, where we have supposed $\pi_v$ is in the complementary series for the nilpotent class containing $N$.  But by Proposition \ref{p:weight-pattern} any unramified $\pi_w$ must also be in the same complementary series.  Moreover, let $s_w$ be the Satake parameter of $\pi_w$  for any unramified $w$.  If $\Re(s_v)$ belongs to the central point (resp. to one of the extraneous classes) then Proposition \ref{p:weight-pattern} implies that $\Re(s_w)$ belongs to the central point (resp. to the same extraneous class) for any unramified $w$.  This completes the proof.

\subsection{Eliminating extraneous parameters} 
Theorem \ref{mainthmArthur} does not imply that the extraneous parameters cannot arise as the local parameters of  a cuspidal (or discrete) automorphic representation $\pi$.  On the contrary, the statement of the theorem asserts that if one local parameter of $\pi$ is extraneous for a given nilpotent orbit $\CO^\vee$, then almost all parameters are extraneous for the same $\CO^{\vee}$.  This would certainly be incompatible with the Arthur Conjectures.  

However, with an additional hypothesis on the local component at $u$ we can rule out extraneous points in the complementary series at all unramified places.
Suppose $\pi_u = AZ(\sigma)$ for an irreducible discrete series representation $\sigma$, with Weil-Deligne parameter $(\varphi_\sigma,N_\sigma)$.  Proposition \ref{AZcusp} asserts that the semisimple parameter is $\CL_{F_u}(\pi_u) = \CL_{F_u}(\sigma)$.  In particular, the weight patterns for $\Frob_u$ for  $\pi_u$ and $\sigma$ are equal.   We apply this in the following Theorem.

\begin{thm}\label{AZs}  Hypotheses are as in Theorem \ref{mainthmArthur}.  Assume moreover that the local component $\pi_u$ is of the form $AZ(\sigma)$ for an irreducible discrete series representation $\sigma$, with Weil-Deligne parameter $(\varphi_\sigma,N)$, where $N$ determines the complementary series for the chosen unramified place $v$.  Then at every split unramified place $w$, the real part of the Satake parameter of $\pi_w$ equals the central point for the complementary series for the nilpotent conjugacy class of $N$.  In particular, no unramified local parameter of $\pi$ is an extraneous point.
\end{thm}

\begin{proof}  We are still assuming that the group $G$ is split at the place $w$. The Weil-Deligne parameter at $u$ gives rise to a simplified Langlands parameter
\[\varphi_{\sigma}^0: W_{K_u}\times SL(2)\to G^\vee
\]
such that $\varphi_\sigma(\Frob_u)=\varphi_{\sigma}^0(\Frob_u) q_u^{h_\sigma^\vee/2}$, where $h_\sigma^\vee$ is a neutral element for a Lie triple of $N$ such that $q_u^{h^\vee/2}=\varphi_{\sigma}^0(\left(\begin{matrix} q_u^{1/2} &0\\0& q_u^{-1/2}\end{matrix}\right))$. Moreover, since $\sigma$ is in particular tempered, $\varphi_{\sigma}^0(\Frob_u)\in G^\vee$ is bounded. Without loss of generality we may assume $\varphi_\sigma(\Frob_u)\in T^\vee$, so that $\varphi_{\sigma}^0(\Frob_u)\in T^\vee_c$ and therefore
\[\Re(\varphi_\sigma(\Frob_u))=q_u^{h_\sigma^\vee/2}.
\]

For any representation $\tau$ of $\LGr$, it follows from the hypotheses of Theorem \ref{mainthmArthur} that the Frobenius weights occurring in the semisimplification of $\CL(\pi,\tau)$  are all integers.  We let $\tau \mapsto \bn(i)(\tau)$ denote the function on representations of $\LGr$ whose value at $\tau$ is $\dim  gr_i(\CL(\pi,\tau))$, for any $i \in \ZZ$.  Thus for any $w$ as in the statement of this theorem, $\bn(i)(\tau) = \bn(i)_w(\tau)$ is the dimension of the weight $i$ (generalized) eigenspace for $\Frob_w$.   

On the other hand, for each $i, k \in \ZZ$ and for any choice of Frobenius element $\Frob_u$ at $u$, we let let $ \tau \mapsto \ba(i,k)(\tau)$ be the function that maps the representation $\tau$ of $\LGr$ to the dimension of the weight $k$ generalized eigenspace of $\Frob_u$ in the fiber at $u$ of $gr_i(\CL(\pi,\tau)$.  There
 is a local monodromy operator $N_i$ at $u$, so that the weight $i-j$ and $i+j$ eigenspaces spaces for $\Frob_u$ in $gr_i(\CL(\pi,\tau))$ are of equal dimension 
 \begin{equation}\label{e:symm}
 \ba(i,i+j)(\tau) = \ba(i,i-j)(\tau).
 \end{equation}
   Indeed, by the purity of the monodromy weight filtration, $N_i^j$ is an isomorphism from one to the other.   We also let $\bn(i)_u(\tau)$ to be the dimension of the weight $i$ generalized eigenspace for $\Frob_u$ in the fiber at $u$ of $\CL(\pi,\tau)_u$.  As we recalled above, this is given by the weight pattern for the discrete series representation $\sigma$.

\smallskip

 We then have
\begin{equation}\label{nijv}
\bn_v(i) = \sum_{k \in \ZZ} \ba(i,k)
\end{equation}
\begin{equation}\label{niju}
\bn_u(i) = \sum_{k \in \ZZ} \ba(k,i)
\end{equation}
as functions of $\tau$.

We claim that (\ref{nijv}) and (\ref{niju}) together with (\ref{e:symm}) imply 

\smallskip

\noindent {\it Property A.} For a representation $\tau$, either 
\begin{enumerate}
\item[(a)]$\bn_u(i)(\tau)= \ba(i,i)(\tau)=\bn_v(i)(\tau)$ for all $i$, or 
\item[(b)]there exists $i_0$ such that $\bn_u(i)(\tau)= \bn_v(i)(\tau)$ for all $i>i_0$ and $\bn_u(i_0)(\tau)> \bn_v(i_0)(\tau)$.
\end{enumerate}

\smallskip

Let us prove this by contradiction. Suppose there exists a representation $\tau$ that violates Property A. Denote for simplicity $a(k,i)=\ba(k,i)(\tau)$, $n_u(i)=\bn_u(i)(\tau)$ but $n_v(i)=\bn_v(i)(\tau)$. There must exist an $i_0$ such that $n_u(i)= n_v(i)$ for all $i>i_0$, but $n_u(i_0)<n_v(i_0)$. Since $\dim\tau$ is finite, there exists $i$ sufficiently large so that $a(i,j)=a(j,i)=0$ for all $j$ (in particular, $n_u(i)=n_v(i)=0$).

By induction suppose that we have proven that $n_v(i)=a(i,i)=n_u(i)$ for all $i> i_1$. This also means that  $a(i,k)=a(k,i)=0$ for all $i>i_1$ and $k\neq i$. In particular, $a(i_1, i')=0$ for all $i'>i_1$. Using the symmetry relation (\ref{e:symm}), we also must have $a(i_1,i'')=0$ for all $i''<i_1$. Hence $a(i_1,j)=0$ for all $j\neq i_1$, which  implies that  \[n_v(i_1)=a(i_1,i_1)\le n_u(i_1).\] 
By our assumption, it follows that $i_1>i_0$ and so $n_v(i_1)=a(i_1,i_1)=n_u(i_1)$. Hence also $a(j,i_1)=0$ for all $j\neq i_1$.

Eventually, we get to $i_1=i_0$, which gives  $n_v(i_0)=a(i_0,i_0)\le n_u(i_0)$. This is a contradiction with the assumption that $n_u(i_0)<n_v(i_0)$. Property A is proven.

\smallskip

In our case, the assumption is that $\pi_v$ belongs to the complementary series of $N$. By Theorem \ref{t:central}, either $\Re(\varphi_{\pi_v}(\Frob_v))={h_\sigma^\vee/2}$ or $\Re(\varphi_{\pi_v}(\Frob_v))$ is one of the extraneous points attached to $N$. It turns out that in the latter case Property A is violated, namely there exists a representation $\tau$ and a weight $i_0$ such that  $\bn_u(i)(\tau)= \bn_v(i)(\tau)$ for all $i>i_0$ but $\bn_u(i_0)(\tau)< \bn_v(i_0)(\tau)$. We give the explicit values in each case, in the notation of Theorem \ref{t:central}.

\begin{enumerate}
\item $G^\vee=Sp(2n)$: $\tau$ is the standard representation, and $i_0=\max\{j\text{ even}| r_j\ge 3\text{ and }\epsilon_j=1\}$. Then $n_v(i_0)=n_u(i_0)+1$.
\item $G^\vee=Spin(n)$: $\tau$ is the standard representation, and $i_0=\max\{i\text{ odd}| r_i\ge 3\text{ and }\epsilon_i=1\}$. Then $n_v(i_0)=n_u(i_0)+1$.
\item $G^\vee=F_4$: $\tau$ is the adjoint representation. 
\begin{enumerate}
\item $\CO^\vee=B_3$: $n_v(i)=n_u(i)=0$ for $i>10$, $n_v(10)=n_u(10)=1$, $n_v(9)=n_v(9)=0$, $n_v(8)=2>1=n_u(8)$.
\item $\CO^\vee=A_1+\wti A_1$:  $n_v(i)=n_u(i)=0$ for $i>4$, $n_v(4)=1>0=n_u(4)$.
\end{enumerate}
\item $G^\vee=E_7$:  $\tau$ is the adjoint representation, $\CO^\vee=D_5(a_1)+A_1$. Then  $n_v(i)=n_u(i)=0$ for $i>10$, $n_v(10)=n_u(10)=1$, $n_v(9)=1>0=n_v(9)$.
\item $G^\vee=E_8$: $\tau$ is the adjoint representation. 
\begin{enumerate}
\item $\CO^\vee=D_6$: $n_v(i)=n_u(i)=0$ for $i>18$, $n_v(18)=n_u(18)=1$, $n_v(17)=n_u(17)=0$, $n_v(16)=2>1=n_v(16)$.
\item $\CO^\vee=A_6$: $n_v(i)=n_u(i)=0$ for $i>13$, $n_v(13)=1>0=n_u(13)$.
\item $\CO^\vee=A_4+A_2$: $n_v(i)=n_u(i)=0$ for $i>9$, $n_v(9)=1>0=n_u(9)$.
\item $\CO^\vee=A_2+2A_1$: $n_v(i)=n_u(i)=0$ for $i>5$, $n_v(5)=1>0=n_u(5)$.
\end{enumerate}

\end{enumerate}

\smallskip

This means that $\varphi_\sigma(\Frob_u)$ and $\varphi_{\pi_v}(\Frob_v)$ must have the same weight pattern, that of $h_\sigma^\vee/2$. As in the previous proof, it follows that the weight pattern at a split unramified place $w$ is the same as that for  $h_\sigma^\vee/2$. Hence by Proposition \ref{p:weight-pattern}, 
\[\Re(\varphi_{\pi_w}(\Frob_w))={h_\sigma^\vee/2}\quad (\text{up to conjugacy}).
\]
\end{proof}

\begin{remark}\label{N0}  It follows from the proof that the monodromy operators $N_i$ are all equal to zero; in other words, $\ba(i,i+j)(\tau) = \ba(i,i-j)(\tau) = 0$ unless $j = 0$.  In particular the local representation at $u$ attached to $\pi$ such that $\pi_u \isoarrow AZ(\sigma)$ has monodromy operator equal to $0$.
On the other hand, the Arthur $SL(2)$ for the discrete series representation $\sigma$  acts trivially.  Thus Theorem \ref{AZs} is compatible with the expectation that the Aubert-Zelevinsky involution exchanges the Deligne and Arthur $SL(2)$ factors.
\end{remark}

\begin{remark}
Notice that using only (\ref{e:symm}), (\ref{nijv}), (\ref{niju}),  it is impossible to conclude purely combinatorially that only case (a) in Property A can occur. Suppose we have a two-dimensional space and the matrix $a(i,j)$ is
\[ a(0,1)=a(0,-1)=1, \quad a(i,j)=0, \text{ for all other }(i,j).
\]
In this example, 
\[\begin{tabular}{|l|c|c|c|}
\hline
$i$ &1&0&-1\\
\hline\hline
$n_u(i)$ (column-sum) &1 &0&1\\
\hline
$n_v(i)$ (row-sum) &0&2&0
\\
\hline
\end{tabular}
\]
\end{remark}

\section*{Appendix: the weight filtration attached to an $SL(2)$-triple}

We give an explicit proof of the following lemma, that may be of independent interest.

\begin{lemma}  Let $(\ve, \vf, \vh)$, $(\ve', \vf', \vh')$ be two $SL(2)$-triples in a simple Lie $C$-algebra $\vg$, and let $(\rho,V) \mapsto (V,Fil_\bullet(V))$, $(\rho,V) \mapsto (V,Fil'_\bullet(V))$ be the corresponding weight filtrations on $Rep(\vg)$.  Suppose for all $(\rho,V)$, and for all $i \in \ZZ$, 
$$\dim Fil_i(V) = \dim Fil'_i(V).$$
Then $(\ve, \vf, \vh)$, $(\ve', \vf', \vh')$ are conjugate under the adjoint action of the semisimple group $\check{G}$ with Lie algebra $\vg$.
\end{lemma}

\begin{proof} Let $\check{\mathfrak a}$ and $\check{\mathfrak a}'$ be the $C$-span in $\vg$ of  $(\ve, \vf, \vh)$ and $(\ve', \vf', \vh')$, respectively. Then $\check{\mathfrak a}\cong\mathfrak{sl}(2,C)$. By $\mathfrak{sl}(2)$ theory, the finite dimensional $\check{\mathfrak a}$-representation $V$ is uniquely determined by $Fil_\bullet (V)$.

If $\vg$ is of classical type, i.e., $\mathfrak{sl}(N,C)$, $\mathfrak {so}(N,C)$, or $\mathfrak{sp}(N,C)$, let $C^N$ denote its standard representation. Let $\check{G}^\#$ denote the Lie group $SL(N,C)$, $SO(N,C)$ is $N$ odd or $O(N,C)$ if $N$ even, and $Sp(N,C)$, corresponding to the Lie algebra $\vg$. The $\check{G}^\#$-conjugacy classes of $\mathfrak{sl}(2)$-subalgebras $\check{\mathfrak a}$ of $\vg$ are parameterized, via the Jordan normal form, by partitions $\lambda=(d_1,d_2,\dots,d_\ell)$ of $N$, with certain additional properties when $\vg=\mathfrak{so}(N)$ or $\mathfrak{sp}(N)$, see \cite[Theorems 5.1.1, 5.1.6]{CM93}. By \cite[Observation 5.1.7]{CM93}, the integers $d_i$ are precisely the dimensions of the irreducible summands of the $\check{\mathfrak a}$-module $C^N$. It follows that $\check{\mathfrak a}$ and $\check{\mathfrak a}'$ are $\check{G}^\#$-conjugate. 

It remains to analyse the case when $\vg=\mathfrak{so}(2n)$. In this case, if the partition $\lambda$ is {\it very even}, i.e., all $d_i$ are even and appear with even multiplicity, then the nilpotent $O(2n)$-orbit labelled by $\lambda$ splits into two $SO(2n)$-orbits $\check\CO_I$ and $\check\CO_{II}$.  We use one of the half-spin modules for $V$ to separate the two orbits. In coordinates for the Lie algebra of type $D_n$, identify the dual of the Cartan subalgebra $\vfh$ with $C^n$ with coordinates $(\ep_1,\dots,\ep_n)$, so that the simple roots are $\pm\ep_1+\ep_2,$ $\ep_2-\ep_3,\dots,\ep_{n-1}-\ep_n$. The two half-spin modules have highest weights $(\pm \frac 12,\frac 12,\dots,\frac 12)$. Choose $S^+$ to be the half-spin module with highest weight $(\frac 12,\frac 12,\dots,\frac 12)$. Its weights are the $n$-tuples with entries $\pm \frac 12$ such that there is an even number of negative entries. Let $\vh_I$ and $\vh_{II}$ be two neutral elements for $\check\CO_I$ and $\check\CO_{II}$, chosen to lie in $\vfh$ and dominant. In coordinates, $\vh_I=(a_1,a_2,\dots,a_n)$ and $\vh_{II}=(-a_1,a_2,\dots,a_n)$, for some odd integers $0<a_1\le a_2\le\dots\le a_n$, see \cite[Recipe 5.2.6]{CM93}. Then the maximal $i$ such that $\dim Fil_i^I(S^+)\neq 0$ is $i_I=\frac 12(a_1+a_2+\dots+a_n)$, while the maximal $i$ such that $\dim Fil_i^{II}(S^+)\neq 0$ is $i_{II}=\frac 12(-a_1+a_2+\dots+a_n)$. This completes the proof for $\mathfrak{so}(2n)$.

\smallskip

For exceptional Lie algebras, using the classification of nilpotent orbits and an explicit computation in {\it Mathematica} reveals that if $V=\vg$ is the adjoint representation, then $\dim Fil_\bullet(\vg)$ uniquely determines the nilpotent orbit, hence the $\mathfrak{sl}(2)$-subalgebra $\check{\mathfrak a}$ up to $\check G$-conjugation. The explicit dimensions for each nilpotent orbit in $E_8$, $F_4$, $G_2$ are listed in the Tables \ref{table:E8-fil}, \ref{table:F4-fil}, \ref{table:G2-fil}. Notice that if the claim holds for $E_8$, it has to hold for all its Levi subalgebras as well, in particular, for $E_6$ and $E_7$, hence the corresponding tables are not presented for economy.

\begin{small}
\begin{longtable}[h]{|c|c|c|}
\hline
Nilpotent &$i_{\max}$ &$\dim Fil_i(\vg)$, $i\ge 0$\\
\hline
$E_8$   &$58$ &$8,0,8,0,7,0,7,0,7,0,7,0,7,0,7,0,6,0,6,0,6,0,6,0,5,0,5,0,4,0,4,0,$\\
&&$4,0,4,0,3,0,3,0,2,0,2,0,2,0,2,0,1,0,1,0,1,0,1,0,1,0,1$\\
\hline
$E_8(a_1)$  &$46$ &$10,0,10,0,9,0,9,0,9,0,9,0,8,0,8,0,7,0,7,0,6,0,6,0,5,0,5,0,4,0,$\\&&$3,0,3,0,3,0,2,0,2,0,1,0,1,0,1,0,1$\\
\hline
$E_8(a_2)$  &$38$ &$12,0,12,0,11,0,11,0,10,0,10,0,9,0,9,0,8,0,7,0,6,0,6,0,4,0,4,0,$\\&&$3,0,2,0,2,0,2,0,1,0,1$\\
\hline
$E_8(a_3)$ &$34$ &$ 14,0,14,0,12,0,12,0,12,0,11,0,9,0,9,0,8,0,7,0,5,0,5,0,4,0,4,0,2,0,1,0,1,0,1$\\
\hline
$E_8(a_4)$ &$28$  &$16,0,16,0,15,0,14,0,13,0,12,0,10,0,10,0,7,0,6,0,4,0,4,0,2,0,2,0,1$\\
\hline
$E_7$  &$34$ &$10,6,7,6,6,6,6,6,6,6,6,4,5,4,5,4,4,4,4,2,3,2,3,2,2,2,2,2,1,0,1,0,1,0,1$\\
\hline
$E_8(b_4)$ &$26$  &$18,0,18,0,16,0,15,0,13,0,13,0,10,0,9,0,7,0,5,0,4,0,3,0,1,0,1$\\
\hline
$E_8(a_5)$  &$22$ &$20,0,20,0,17,0,16,0,15,0,14,0,10,0,8,0,5,0,4,0,3,0,2$\\
\hline
$E_7(a_1)$  &$26$ &$12,8,9,8,8,8,8,6,7,6,7,6,5,4,5,4,4,2,3,2,2,2,2,0,1,0,1$\\
\hline
$E_8(b_5)$  &$22$ &$22,0,22,0,18,0,17,0,15,0,12,0,9,0,9,0,6,0,3,0,1,0,1$\\
\hline
$D_7$ &$22$  &$12,10,9,10,8,8,8,8,7,8,7,6,6,4,3,4,2,2,2,2,1,2,1$\\
\hline
$E_8(a_6)$  &$18$ &$24,0,24,0,21,0,20,0,15,0,12,0,9,0,6,0,3,0,2$\\
\hline
$E_7(a_2)$  &$22$  &$14,10,11,10,9,8,9,8,8,6,7,4,5,4,5,4,3,2,2,0,1,0,1$\\
\hline
$E_6A_1$  &$22$ &$14,12,11,8,9,8,9,8,9,6,6,4,5,4,5,4,4,2,1,0,1,0,1$\\
\hline
$D_7(a_1)$  &$18$ &$26,0,25,0,21,0,19,0,16,0,13,0,7,0,6,0,3,0,1$\\
\hline
$E_8(b_6)$  &$16$ &$28,0,28,0,24,0,20,0,15,0,12,0,6,0,4,0,1$\\
\hline
$E_7(a_3)$  &$18$ &$16,12,13,10,11,10,10,8,8,8,7,4,4,2,4,2,2,0,1$\\
\hline
$E_6(a_1)A_1$  &$16$ &$16,14,15,12,13,10,10,8,9,6,6,4,4,2,2,0,1$\\
\hline
$A_7$  &$15$ &$16,14,13,14,12,12,9,10,8,6,5,4,4,2,1,2$\\
\hline
$D_7(a_2)$ &$14$   &$16,16,15,14,13,12,10,10,7,6,4,4,2,2,1$\\
\hline
$E_6$  &$22$ &$32,0,18,0,17,0,17,0,17,0,10,0,9,0,9,0,8,0,1,0,1,0,1$\\
\hline
$D_6$ &$18$   &$20,12,10,12,9,12,9,8,8,8,8,4,2,4,2,4,1,0,1$\\
\hline
$D_5A_2$ &$14$  &$34,0,33,0,25,0,20,0,15,0,10,0,3,0,1$\\
\hline
$E_6(a_1)$ &$16$  &$34,0,26,0,25,0,18,0,17,0,10,0,8,0,2,0,1$\\
\hline
$E_7(a_4)$  &$14$ &$20,16,17,14,13,10,11,8,8,6,6,2,2,0,1$\\
\hline
$A_6A_1$  &$12$  &$ 20,16,17,14,15,12,12,8,7,4,4,2,3$\\
\hline
$D_6(a_1)$   &$14$  &$22,16,16,16,10,12,10,8,8,8,3,4,1,0,1$\\
\hline
$A_6$    &$12$   &$38,0,32,0,27,0,24,0,11,0,8,0,3$\\
\hline
$E_8(a_7)$  &$10$ &$40,0,40,0,30,0,20,0,10,0,4$\\
\hline
$D_5A_1$  &$14$   &$22,18,16,12,14,10,10,8,9,6,6,2,1,0,1$\\
\hline
$E_7(a_5)$ &$10$  &$24,18,21,18,15,12,11,6,6,2,3$\\
\hline
$E_6(a_3)A_1$  &$10$ &$24,20,21,16,17,12,10,6,6,2,2$\\
\hline
$D_6(a_2)$ &$10$  &$24,20,18,20,15,12,10,8,3,4,2$\\
\hline
$D_5(a_1)A_2$  &$10$ &$24,22,21,18,14,14,8,8,4,2,1$\\
\hline
$A_5A_1$   &$10$  &$24,22,18,18,16,14,9,6,4,4,1$\\
\hline
$A_4A_3$ &$9$    &$24,24,21,20,15,12,9,6,3,2$\\
\hline
$D_5$    &$14$   &$48,0,27,0,26,0,18,0,17,0,10,0,1,0,1$\\
\hline
$E_6(a_3)$  &$10$ &$50,0,36,0,33,0,18,0,10,0,2$\\
\hline
$D_4A_2$   &$10$  &$50,0,42,0,28,0,21,0,7,0,1$\\
\hline
$A_4A_2A_1$  &$8$ &$28,24,25,18,18,10,8,4,3$\\
\hline
$D_5(a_1)A_1$  &$10$ &$28,24,22,14,15,12,12,6,4,0,1$\\
\hline
$A_5$    &$10$   &$34,18,17,18,16,16,9,2,8,2,1$\\
\hline
$A_4A_2$   &$8$  &$54,0,48,0,30,0,16,0,3$\\
\hline
$A_42A_1$ &$8$   &$28,28,24,20,15,12,6,4,1$\\
\hline
$D_5(a_1)$  &$10$ &$34,24,19,16,11,16,10,8,2,0,1$\\
\hline
$2A_3$   &$7$   &$32,28,22,24,16,8,6,4$\\
\hline
$A_4A_1$  &$8$   &$34,26,25,18,17,10,8,2,1$\\
\hline
$D_4(a_1)A_2$  &$6$ &$64,0,56,0,28,0,8$\\
\hline
$D_4A_1$  &$10$ &$38,26,17,12,15,12,15,6,1,0,1$\\
\hline
$A_3A_2A_1$ &$6$  &$36,30,30,20,15,6,5$\\
\hline
$A_4$  &$8$   &$68,0,44,0,33,0,12,0,1$\\
\hline
$A_3A_2$ &$6$  &$38,32,27,24,11,8,3$\\
\hline
$D_4(a_1)A_1$  &$6$ &$40,32,31,18,15,6,2$\\
\hline
$A_32A_1$  &$6$ &$40,36,27,22,12,6,1$\\
\hline
$2A_22A_1$ &$5$  &$40,40,30,20,10,4$\\
\hline
$D_4$   &$10$   &$80,0,28,0,27,0,27,0,1,0,1$\\
\hline
$D_4(a_1)$ &$6$   &$82,0,54,0,27,0,2$\\
\hline
$A_3A_1$   &$6$ &$50,34,26,20,16,2,1$\\
\hline
$2A_2A_1$  &$5$ &$50,36,33,18,10,2$\\
\hline
$2A_2$   &$4$  &$92,0,64,0,14$\\
\hline
$A_23A_1$ &$4$   &$52,42,35,14,7$\\
\hline
$A_3$   &$6$  & $68,32,13,32,12,0,1$\\
\hline
$A_22A_1$ &$4$  & $54,48,30,16,3$\\
\hline
$A_2A_1$   &$4$ &$68,44,33,12,1$\\
\hline
$4A_1$   &$3$  &$64,56,28,8$\\
\hline
$A_2$  &$4$    &$134,0,56,0,1$\\
\hline
$3A_1$   &$3$  &$82,54,27,2$\\
\hline
$2A_1$ &$2$    &$92,64,14$\\
\hline
$A_1$   &$2$   &$134,56,1$\\
\hline
$1$     &$0$ &$248$\\
\hline

\caption{Filtrations for the nilpotent orbits of $E_8$.}\label{table:E8-fil}
\end{longtable}
\end{small}

\begin{small}
\begin{longtable}[h]{|c|c|c|}
\hline
Nilpotent &$i_{\max}$ &$\dim Fil_i(\vg)$, $i\ge 0$\\
\hline

$F_4$ &$22$ &$4,0,4,0,3,0,3,0,3,0,3,0,2,0,2,0,1,0,1,0,1,0,1$\\
\hline

$F_4(a_1)$ &$14$ &$6,0,6,0,5,0,4,0,3,0,3,0,1,0,1$\\
\hline

$F_4(a_2)$ &$10$ &$8,0,8,0,5,0,4,0,3,0,2$\\
\hline

$C_3$ &$10$ &$6,4,3,4,2,2,2,2,1,2,1$\\
\hline

$B_3$ &$10$ &$10,0,7,0,6,0,6,0,1,0,1$\\
\hline

$F_4(a_3)$ &$6$  &$12,0,12,0,6,0,2$\\
\hline

$C_3(a_1)$ &$6$  &$8,6,5,6,2,2,1$\\
\hline

$A_1+\widetilde A_2$ &$5$   &$8,8,5,4,3,2$\\
\hline

$B_2$ &$6$ &$12,4,6,4,5,0,1$\\
\hline

$\widetilde A_1+A_2$ &$4$ &$12,6,9,2,3$\\
\hline

$\widetilde A_2$ &$4$  &$22,0,8,0,7$\\
\hline

$A_2$  &$4$ &$22,0,14,0,1$\\
\hline

$A_1+\widetilde A_1$ &$3$ &$12,12,6,2$\\
\hline

$\widetilde A_1$ &$2$ &$22,8,7$\\
\hline

$A_1$&$2$ &$22,14,1$\\
\hline

$1$ &$0$ &$52$\\
\hline

\caption{Filtrations for the nilpotent orbits of $F_4$.}\label{table:F4-fil}
\end{longtable}
\end{small}

\begin{small}
\begin{longtable}[h]{|c|c|c|}
\hline
Nilpotent &$i_{\max}$ &$\dim Fil_i(\vg)$, $i\ge 0$\\
\hline

$G_2$ &$10$ &$2,0,2,0,1,0,1,0,1,0,1$\\
\hline
$G_2(a_1)$ &$4$ &$4,0,4,0,1$ \\
\hline

$\widetilde A_1$ &$3$ &$4,2,1,2$\\
\hline

$A_1$ &$2$ &$4,4,1$\\
\hline

$1$ &$0$ &$14$\\
\hline

\caption{Filtrations for the nilpotent orbits of $G_2$.}\label{table:G2-fil}
\end{longtable}
\end{small}

\end{proof}

\end{document}